\documentclass[3p,preprint,sort&compress,lefttitle]{elsarticle} 
\usepackage{xcolor}

\usepackage[bookmarksnumbered=true,colorlinks=true,citecolor=green!50!black,linkcolor=red!80!black]{hyperref}

\usepackage{amsmath,amssymb,amsthm}

\usepackage[mathlines]{lineno}

\newcommand*\patchAmsMathEnvironmentForLineno[1]{
	\expandafter\let\csname old#1\expandafter\endcsname\csname #1\endcsname
	\expandafter\let\csname oldend#1\expandafter\endcsname\csname end#1\endcsname
	\renewenvironment{#1}
	{\linenomath\csname old#1\endcsname}
	{\csname oldend#1\endcsname\endlinenomath}}
\newcommand*\patchBothAmsMathEnvironmentsForLineno[1]{
	\patchAmsMathEnvironmentForLineno{#1}
	\patchAmsMathEnvironmentForLineno{#1*}}
\patchBothAmsMathEnvironmentsForLineno{equation}
\patchBothAmsMathEnvironmentsForLineno{align}
\patchBothAmsMathEnvironmentsForLineno{flalign}
\patchBothAmsMathEnvironmentsForLineno{alignat}
\patchBothAmsMathEnvironmentsForLineno{gather}
\patchBothAmsMathEnvironmentsForLineno{multline}

\usepackage[todonotes={textsize=tiny}]{changes}
\definechangesauthor[name={Zhicheng HU}, color=red]{zc}

\setdeletedmarkup{{\color{blue} \sout{#1}}}

\usepackage{subfig}

\usepackage{cleveref}
\crefname{equation}{}{}
\crefname{algorithm}{Algorithm}{Algorithm}
\crefname{table}{Table}{Table}
\crefname{figure}{\figurename}{\figurename}

\usepackage{booktabs}  

\graphicspath{{images}}

\newcommand\pd[2]{\dfrac{\partial {#1}}{\partial {#2}}}
\newcommand\pdtw[2]{\dfrac{\partial^2 {#1}}{\partial {#2}^2}}

\newcommand\dd{\,\mathrm{d}}

\newcommand\mL{\mathcal{L}}

\begin{document}
	
	\hypersetup{citecolor=green!50!black,linkcolor=red!80!black}
	\begin{frontmatter}
		
		\title{Thermal analysis of dual-phase-lag model in a two-dimensional plate subjected to a heat source moving along elliptical trajectories}
		
		\author[]{Kaiyuan Chen}
		\ead{chenkaiyuan2000@163.com}
		
		\author[]{Zhicheng Hu\corref{cor1}}
		\ead{huzhicheng@nuaa.edu.cn}
		
		\cortext[cor1]{Corresponding author.}
		\affiliation[]{organization={School of Mathematics, Nanjing University
				of Aeronautics and Astronautics},
			city={Nanjing},
			postcode={211106}, 
			country={China}}
		
		\affiliation[]{organization={Key Laboratory of Mathematical
				Modelling and High Performance Computing of Air Vehicles
				(NUAA), MIIT},
			city={Nanjing},
			postcode={211106}, 
			country={China}}

		\begin{abstract}
			In this paper, we focus on the study of heat transfer behavior for
			the dual-phase-lag heat conduction model, which describes the
			evolution of temperature in a two-dimensional rectangular plate
			caused by the activity of a point heat source moving along
			elliptical trajectories. At first, Green's function approach is
			applied to derive the analytical solution of temperature for the
			given model. Based on the series representation of this analytical
			solution, the thermal responses for the underlying heat transfer
			problem, including the relations between the moving heat source and
			the concomitant temperature peak, the influences of the pair of
			phase lags and the angular velocity of heat source on temperature,
			are then investigated, analyzed and discussed in detail for three
			different movement trajectories. Compared with the results revealed
			for the common situation that the heat source moves in a straight
			line with a constant speed, the present results show quite
			distinctive thermal behaviors for all cases, which subsequently can
			help us to better understand the internal mechanism of the
			dual-phase-lag heat transfer subjected to a moving heat source with
			curved trajectory.
		\end{abstract}

		
	\end{frontmatter}

\section{Introduction}
\label{sec:intro}

Heat conduction is a natural phenomenon that can be experienced so
commonly in our daily life.
With the rapid development of localized high-intensity heat source
technology such as laser beam and the wide application of moving heat
source in numerous fields, including cutting, welding, drilling, laser
hardening/forming, plasma spraying, and even modern medicine such as
laser surgery and hyperthermia treatment, the heat conduction process
in a finite medium involving a moving heat source has been of
particular interest by engineers and scientists in the past few
decades, see e.g. \cite{hahn2012inbook, panas2014inbook,
	champagne2020numerical, 2022, ma2019thermal, kabiri2021thermal,
	yang2023transient, partovi2023analytical} and references therein.
In order to explore thermal properties of the medium caused by the
activity of heat source, one of the core tasks is to precisely and
efficiently determine the responsive temperature in the medium,
as it is a fundamental physical quantity and many other thermal
quantities of interest, such as heat flux, thermal stress, thermal
deformation, and biological tissue damage, could then be evaluated
based on it \cite{mirkoohi2018thermal, bian2019finite, 2018,
	ma2019thermal, kabiri2021thermal}. Lots of mathematical models and
methods have been established to predict the temperature field for
these moving heat source problems.

Usually, the classical Fourier's law, assuming that the heat flux is
proportional to the temperature gradient, can well describe the heat
conduction process. It follows that the evolution of temperature in
the medium heated by a moving heat source would satisfy the classical
heat conduction equation with a time-dependent localized source term,
which is introduced to simulate the contribution of the moving heat
source.
Accordingly, many heat transfer studies have been carried out based on
these classical heat conduction models. Several examples are given
hereunder. In the work of Araya and Gutierrez
\cite{araya2006analytical}, an analytical solution of the transient
temperature distribution in a three-dimensional (3D) finite domain
subjected to a moving laser beam of two types was
presented. Considering applications in laser material processing, Van
Elsen et al. \cite{van2007solutions} described both analytical and
numerical solutions of the heat conduction equation with a localized
moving heat source of any type. A comprehensive investigation for the
geometrical effects of a moving heat source on the temperature
distribution was carried out by Akbari et
al. \cite{akbari2011geometrical}. Several heat source models including
surface and volumetric sources were discussed by Winczek
\cite{winczek2016influence} using a submerged arc welding problem. To
investigate underwater welding and cutting, Salmi et
al. \cite{salimi2016analytical} derived an analytical solution for the
thermal problems with varying boundary conditions around a moving
source. For numerical simulation of the temperature field due to a
moving heat source in arc welding, Champagne and Pham
\cite{champagne2020numerical} developed a meshless element-free
Galerkin method which was compared with the finite element method as
well as experimental data. To improve accuracy without loss of
efficiency, a dynamic node reconfiguration strategy for meshfree 
methods was proposed by Khosravifard et
al. \cite{khosravifard2020meshfree}.
Some adaptive moving mesh methods proposed in \cite{hu2018numerical,
	hu2020heat} by Hu and his co-worker for the heat conduction problem
with multiple moving heat sources can improve numerical efficiency
greatly too.

Nevertheless, it has also been observed in some experiments that the
classical Fourier's law might be inadequate to describe the heat
conduction process accurately, especially when the temporal or spatial
characteristic size reduces to microscopic scale
\cite{wang2007heat}. To remedy such a deficiency, Cattaneo
\cite{cattaneo1958form} and Vernotte \cite{vernotte1958paradoxes}
independently proposed a hyperbolic heat conduction model, now known
as the CV model or the thermal wave model, by introducing a hysteretic
time for the heat flux such that the CV model could account for the
non-Fourier heat conduction phenomenon and make the heat propagate as
waves with a finite speed. Although it has been validated both
theoretically and experimentally that the CV model might be able to
give better description of the heat propagation than the classical
model in the fields such as ultra-fast lasers, there are still some
shortcomings which shall limit the application of the CV model
\cite{bai1995hyperbolic, tzou2014macro}. For example, the CV model
might be inconsistent with the second law of thermodynamics due to its
hyperbolicity. In order to partially resolve these issues of
microscopic heat conduction, Tzou \cite{tzou} also introduced a time
lag for the temperature gradient, which together with the heat flux
time lag allows the temperature gradient ahead of the heat flux or
vice versa in the transient heat transfer process. The resulting heat
conduction model is referred to as the dual-phase-lag (DPL) model,
which could be consistent with the second law of thermodynamics,
provided that both time lags are non-negative
\cite{kovacs2018thermo}. In recent years, the DPL model has attracted
increasing attention and been widely utilized with various newly
developed models in not only the engineering fields but also the
fields of modern medicine where the non-Fourier effect can not be
negligible \cite{tzou2001temperature, ghazanfarian2012investigation,
	zhou2020dual, li2021fractional, ma2021theoretical, ghasemi2022dual,
	ramos2023mathematical, ma2023thermodynamic}. Particularly, many
researchers have taken the DPL model into consideration for the heat
transfer analysis in a medium, whether material of metal, non-metal or
biological tissue, when it is subjected to a moving heat source, see
e.g. \cite{lee2016numerical, 2018, ma2019thermal, 2022,
	yang2023transient, partovi2023analytical}.

As shown in most studies of moving heat source problems,
e.g. \cite{hahn2012inbook, panas2014inbook, salimi2016analytical,
	2018, 2022}, the heat source is assumed to move along a straight
line with a constant speed. For the situation that the heat source
moves along a curved trajectory, we have tried our best and only
retrieved a few studies involving the classical heat conduction
models. See below for instances.
Applying Green's function method, Kidawa-Kukla \cite{2008}
obtained analytical solutions of temperature in a rectangular plate
heated by a moving heat source with elliptical trajectories.
Again using Green's function method, Elsheikh et
al. \cite{elsheikh2018thermal} investigated thermal effects on the
deflection and stresses in a thin-wall workpiece during a turning
process with an axisymmetric heat source. In the works of Hu and his
co-worker \cite{hu2018numerical, hu2020heat} for efficient numerical
simulation, the moving heat source is also allowed to move with
various curved trajectories.

Motivated by the facts that there are practical applications which
indeed involve the heat source moving in a non-linear trajectory, and
such a moving heat source might give rise to quite different heat
transfer behaviors, we in this paper concentrate on the analysis of
heat transfer generated by the two-dimensional (2D) DPL 
model when the medium is heated by a moving heat source with a curved
trajectory. The point heat source described by the Dirac delta
function with three different trajectories, namely, line segment
trajectory, circular trajectory, and elliptical trajectory, is taken
into account. Following Green's function approach \cite{cole2011book},
the analytical solution of temperature in series form is established.
The convergence of the series solution is validated via the distance
between the positions of the moving heat source and the concomitant
temperature peak. The temperature responses of the underlying DPL
model induced by different trajectories, various pairs of phase lags,
and various angular velocities of heat source are then
investigated. Besides observing similar heat transfer phenomena to the
situation of the heat source moving in a straight line, the present
results also exhibit some quite distinctive thermal behaviors caused
by the curved motion of heat source.

The rest of the paper is arranged as follows. In \cref{sec:model}, the
DPL model in a two-dimensional rectangular plate irradiated by a
moving point heat source is described. The detailed formulation of the
analytical solution for the underlying DPL model is established in
\cref{sec:a-sol}. The analysis and discussion are then carried out in
\cref{sec:result} based on the results calculated from the analytical
solution. Finally, some conclusions are given in the last section.
\section{Dual-phase-lag model with a moving heat source}
\label{sec:model}

\begin{figure}[!htbp] 
	\centering
	\includegraphics[width=0.5\textwidth]{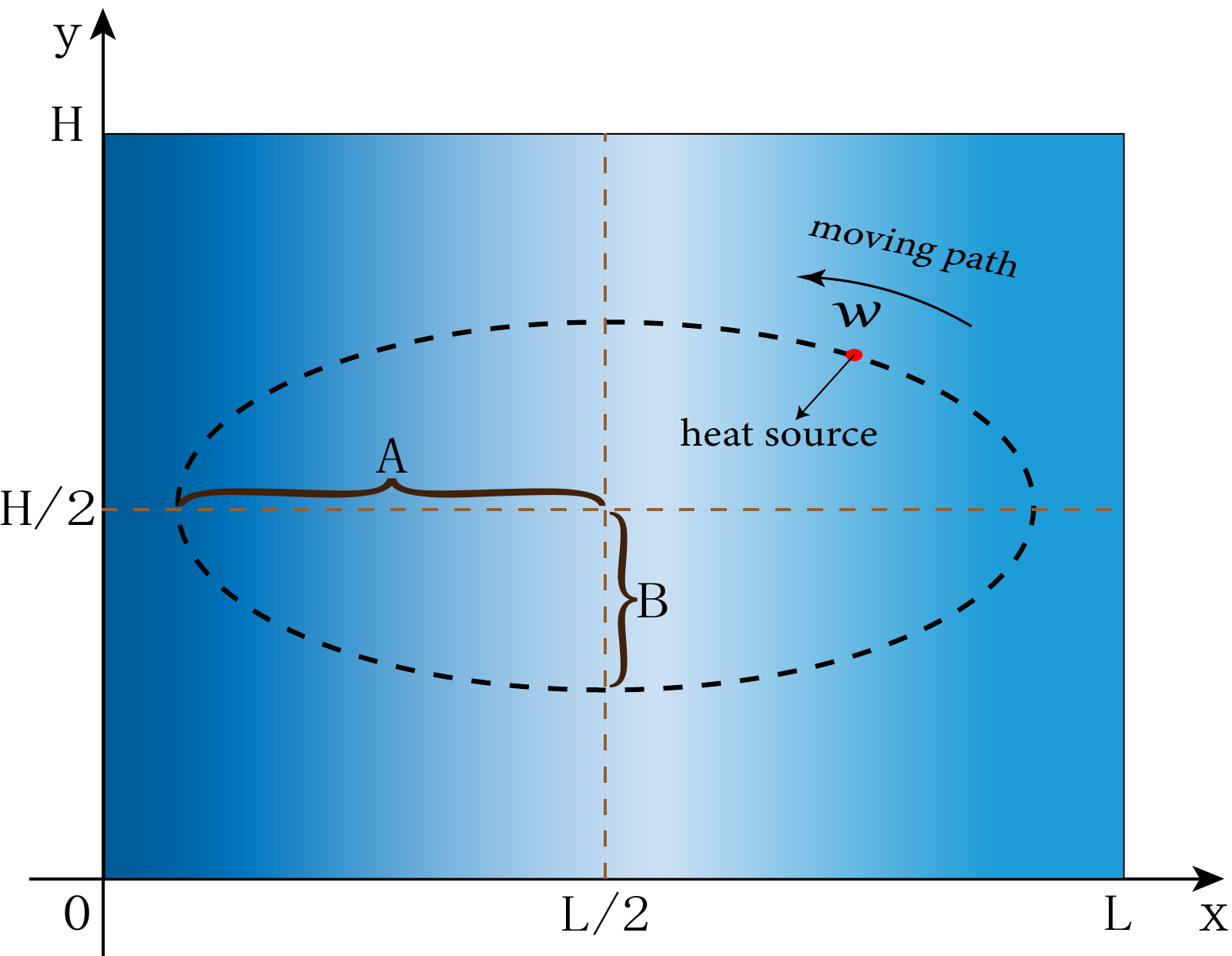} 
	\caption{Schematic diagram of the rectangular plate with a moving point heat source.}
	\label{fig:sketch} 
\end{figure}

Consider a thin rectangular plate that is heated by a moving line heat
source oriented perpendicular to the plate. When the plate is
homogeneous and isotropic, the heat transfer in it can be simplified
to be viewed as a 2D heat transfer problem, in which the line heat
source reduces to a point heat source \cite{hahn2012inbook,
	panas2014inbook}.
For a moving point heat source of constant strength $\Theta$ with the
position given by time-dependent functions $\bar{x}(t)$ and
$\bar{y}(t)$, it can be well modeled by using Dirac delta function as
\begin{equation}
	\label{eq:heat-source}
	Q(x,y,t) = \Theta \delta(x-\bar{x}(t)) \delta(y-\bar{y}(t)).
\end{equation}
As shown in \cref{fig:sketch} for the sketch of the problem, we assume
in this paper that the point heat source moves along an elliptical
trajectory around the center of the plate with a constant angular
velocity $w$. So we have
\begin{equation}
	\label{eq:trajectory}
	\bar{x}(t) = L/2 + A \cos w t, \qquad \bar{y}(t) = H/2 + B \sin wt, 
\end{equation}
where $L, H$, and $A, B$ are length, height of the plate, and lengths
of semi-major and semi-minor axes of the ellipse, respectively.

Let $\boldsymbol{q}(x,y,t)$ and $T(x,y,t)$ be heat flux vector and
temperature, respectively, at position $(x,y)$ of the plate at time
$t$. To get a more accurate prediction of thermal process than the
classical Fourier's law in the situation,
when the fast-transient effects and micro-structural interactions
could not be neglected, Tzou \cite{tzou} introduced two phase lags
$\tau_{q}$ and $\tau_{T}$ into Fourier's law such that the
constitutive relation between $\boldsymbol{q}$ and $T$ becomes
\begin{equation}
	\label{eq:nonFourier-law}
	\boldsymbol{q}(x,y,t+\tau_{q}) = -k \nabla T(x,y,t+\tau_{T}),
\end{equation}
in which $k$ is the thermal conductivity, and
$\nabla= (\frac{\partial}{\partial x}, \frac{\partial}{\partial
	y})^{T}$ is the gradient operator. Applying the Taylor expansion to
both sides of the DPL constitutive relation \cref{eq:nonFourier-law},
a first-order approximation of it can be obtained as follows
\begin{equation}
	\label{eq:1st-dpl-law}
	\boldsymbol{q}(x,y,t) + \tau_{q} \frac{\partial \boldsymbol{q}(x,y,t)}{\partial t} = - k \left[ \nabla T (x,y,t) + \tau_{T} \frac{\partial \nabla T(x,y,t)}{\partial t} \right],
\end{equation}
which is still called the DPL constitutive relation and also known as
the Jeffreys-type constitutive relation \cite{wang2007heat,
	tamma2014inbook}.

Based on the DPL constitutive relation \cref{eq:1st-dpl-law} and the
conservation of energy, the evolution of temperature $T$ in the plate
would satisfy the DPL heat conduction equation \cite{wang2007heat,
	2018, 2022} given by
\begin{equation}
	\label{eq:dpl-heat-conduction}
	\nabla \cdot \nabla T + \tau_{T} \pd{}{t}(\nabla \cdot \nabla T) + \frac{1}{k} \left( Q + \tau_{q} \pd{Q}{t} \right) = \frac{1}{\alpha} \left( \pd{T}{t} + \tau_{q} \pdtw{T}{t} \right), 
\end{equation}
where $\alpha= k/(\rho C)$ is the thermal diffusivity of the plate
with $\rho$ be the material density and $C$ be the heat capacity. In
particular, the above heat conduction equation degenerates to the CV
heat conduction equation when $\tau_{T}=0$, and to the classical heat
conduction equation when $\tau_{q} = \tau_{T}=0$, respectively.

It remains to complement the initial and boundary conditions for the
DPL heat conduction equation \cref{eq:dpl-heat-conduction}. As
investigated in \cite{2022}, we suppose that the environmental
temperature surrounding the plate is $T_{0}$. Then, the initial
condition at time $t=0$ could be set as
\begin{equation}
	\label{eq:initial-condition}
	T(x,y,0) = T_{0}, \qquad 
	\pd{T(x,y,0)}{t} = 0,
\end{equation}
and the boundary condition throughout the time would be adopted as
\begin{equation}
	\label{eq:boundary-condition}
	T(0,y,t) = T(L, y,t) = T(x,0,t) = T(x,H,t) = T_{0}.
\end{equation}

For convenience of discussion, we shall directly view the DPL heat
conduction equation \cref{eq:dpl-heat-conduction}, together with the
initial and boundary conditions
\cref{eq:initial-condition,eq:boundary-condition}, as a dimensionless
model
in the current work. Accordingly, 
we further assume the environmental temperature $T_{0} = 0$ for
simplicity.

\section{Formulation of analytical solution}
\label{sec:a-sol}

This section is devoted to find the analytical solution of the DPL
heat conduction problem
\cref{eq:dpl-heat-conduction,eq:initial-condition,eq:boundary-condition}
with the moving heat source function \cref{eq:heat-source} by using
Green's function method.

\subsection{Integral form of solution in terms of Green's function}
\label{sec:a-sol-green}

Let $\mL$ be the modified heat flux linear operator defined by
\begin{equation}
	\label{eq:operator-L}
	\mL = \frac{1}{\alpha} \left( \pd{}{t} + \tau_{q} \pdtw{}{t} \right) - \left(\pdtw{}{x} + \pdtw{}{y} \right) - \tau_{T} \left( \pd{^{3}}{t\partial x^{2}} + \pd{^{3}}{t\partial y^{2}} \right).
\end{equation}
It follows that the DPL heat conduction equation
\cref{eq:dpl-heat-conduction} can be rewritten as
\begin{equation}
	\label{eq:dpl-heat-eqnew}
	\mL T(x,y,t) = \frac{1}{k} \left( Q + \tau_{q} \pd{Q}{t} \right) \triangleq \bar{Q}(x,y,t).
\end{equation}
Specifically, for the moving heat source described by
\cref{eq:heat-source}, we have
\begin{equation}
	\label{eq:rhs-dpl-heat-eq}
	\bar{Q}(x,y,t) = \frac{\Theta}{k}\left( \delta(x-\bar{x}(t)) \delta(y-\bar{y}(t)) -\tau_{q} \pd{\bar{x}(t)}{t} \delta(1,x-\bar{x}(t))\delta(y-\bar{y}(t)) - \tau_{q} \pd{\bar{y}(t)}{t} \delta(x-\bar{x}(t)) \delta(1,y-\bar{y}(t)) \right),
\end{equation}
where $\delta(1,\cdot)$ represents the derivative of the Dirac delta
function.

Following Green's function approach employed in
\cite{Frankel,2018,2022}, we would first introduce the operator $\mL'$
by replacing variables $x$, $y$, and $t$ in the operator $\mL$ with
$\xi$, $\eta$, and $\tau$, respectively, that is,
\begin{equation}
	\label{eq:operator-Lxi}
	\mL' = \frac{1}{\alpha} \left( \pd{}{\tau} + \tau_{q} \pdtw{}{\tau} \right) - \left(\pdtw{}{\xi} + \pdtw{}{\eta} \right) - \tau_{T} \left( \pd{^{3}}{\tau\partial \xi^{2}} + \pd{^{3}}{\tau\partial \eta^{2}} \right).
\end{equation}
The Green's function applied to solve the considered heat conduction
problem could then be introduced by using the formal adjoint operator
of the above operator, which is defined as
\begin{equation}
	\label{eq:operator-L-adjoint}
	\mL'^{*} = \frac{1}{\alpha} \left( - \pd{}{\tau} + \tau_{q} \pdtw{}{\tau} \right) - \left(\pdtw{}{\xi} + \pdtw{}{\eta} \right) + \tau_{T} \left( \pd{^{3}}{\tau\partial \xi^{2}} + \pd{^{3}}{\tau\partial \eta^{2}} \right).
\end{equation}
Precisely speaking, the Green's function $G(x,y,t|\xi,\eta,\tau)$ is
determined as the solution of the differential equation
\begin{equation}
	\label{eq:green-eq}
	\mL'^{*} G(x,y,t|\xi,\eta,\tau) = \delta(x-\xi) \delta(y-\eta) \delta(t-\tau),
\end{equation}
with homogeneous boundary condition
\begin{equation}
	\label{eq:green-eq-bc}
	G(x,y,t|0,\eta,\tau) = G(x,y,t|L,\eta,\tau) = G(x,y,t|\xi,0,\tau) = G(x,y,t|\xi,H,\tau) = 0,
\end{equation}
and additional condition
\begin{equation}
	\label{eq:green-eq-ic}
	G(x,y,t|\xi,\eta,\tau) = \pd{G(x,y,t|\xi,\eta,\tau)}{\tau} = 0, \quad t<\tau.
\end{equation}

Now let us consider the integral given by
\begin{equation}
	\label{eq:integral-def}
	I = \lim_{\varepsilon\to 0} \int_{0}^{t+\varepsilon} \int_{0}^{L} \int_{0}^{H} G(x,y,t|\xi,\eta,\tau) \mL'T(\xi,\eta,\tau) \dd \eta \dd \xi \dd \tau.
\end{equation}
According to the integration by parts, we could rewrite it as
\begin{equation}
	\label{eq:integral-sum}
	I = \lim_{\varepsilon \to 0} \left(I_{1} + I_{2} + I_{3} + I_{4}\right),
\end{equation}
where
\begin{align}
	\label{eq:integral-I1}
	& I_{1} = \int_{0}^{L} \int_{0}^{H} \left.\left( \frac{GT}{\alpha} + \frac{\tau_{q} G}{\alpha}\pd{T}{\tau} - \frac{\tau_{q}}{\alpha} \pd{G}{\tau} T - \tau_{T} G\left( \pdtw{T}{\xi} +\pdtw{T}{\eta}\right) \right)\right|_{0}^{t+\varepsilon}\dd \eta \dd \xi,\\
	& I_{2} = \int_{0}^{t+\varepsilon}\int_{0}^{H} \left.\left( -G \pd{T}{\xi} + \pd{G}{\xi} T + \tau_{T} \pd{G}{\tau} \pd{T}{\xi} - \tau_{T} \pd{^{2}G}{\xi\partial \tau} T \right)\right|_{0}^{L} \dd \eta \dd \tau, \label{eq:integral-I2}\\
	& I_{3} = \int_{0}^{t+\varepsilon}\int_{0}^{L} \left.\left( -G \pd{T}{\eta} + \pd{G}{\eta} T + \tau_{T} \pd{G}{\tau} \pd{T}{\eta} - \tau_{T} \pd{^{2}G}{\eta\partial \tau} T \right)\right|_{0}^{H} \dd \xi \dd \tau, \label{eq:integral-I3}\\
	& I_{4} = \int_{0}^{t+\varepsilon}\int_{0}^{L}\int_{0}^{H} T(\xi,\eta,\tau) \mL'^{*}G(x,y,t|\xi,\eta,\tau) \dd \eta \dd \xi \dd \tau. \label{eq:integral-I4}
\end{align}
Since the initial condition of the temperature, the additional
condition of the Green's function, and the boundary conditions of them
are all assumed to be zero, it is easy to show that
$I_{1}=I_{2}=I_{3} = 0$. As for the integral $I_{4}$, substituting
\cref{eq:green-eq} into \cref{eq:integral-I4} and using the properties
of the Dirac delta function, we obtain $I_{4} =
T(x,y,t)$. Therefore, we can deduce that
\begin{equation}
	\label{eq:integral-form-sol}
	T(x,y,t) = I = \int_{0}^{t} \int_{0}^{L} \int_{0}^{H} G(x,y,t|\xi,\eta,\tau) \bar{Q}(\xi,\eta,\tau) \dd \eta \dd \xi \dd \tau,
\end{equation}
by noting that $\mL' T(\xi,\eta,\tau) = \bar{Q}(\xi,\eta,\tau)$.

\subsection{Series representation of solution}
\label{sec:a-sol-series}

In order to further represent the temperature into a series from the
integral expression \cref{eq:integral-form-sol}, we shall return to
the problem \cref{eq:green-eq,eq:green-eq-bc,eq:green-eq-ic} to
determine the series of Green's function as a premise.

At first, it is straightforward to verify that, for any positive
integers $m$ and $n$, the trigonometric functions
\begin{equation}
	\label{eq:trig-eig-fun}
	\sin\mu_{m} \xi \quad \text{and} \quad \sin \gamma_{n} \eta,
\end{equation}
with $\mu_{m} = m\pi/L$ and $\gamma_{n}=n\pi/H$, respectively, are
eigenfunctions of the operator $\mL'^{*}$, that satisfy the boundary
condition \cref{eq:green-eq-bc}.
Similar to \cite{2018,2022}, it follows that the Green's function of
the current problem can be expanded into a series in terms of these
eigenfunctions, that is,
\begin{equation}
	\label{eq:green-fun-series-def}
	G(x,y,t|\xi,\eta,\tau) = \sum_{m=1}^{\infty}\sum_{n=1}^{\infty} \frac{4}{LH} \bar{G}_{mn}(x,y,t|\tau) \sin\mu_{m}\xi\sin\gamma_{n}\eta.
\end{equation}
From the orthogonality of eigenfunctions, we get
\begin{equation}
	\label{eq:green-fun-series-coe}
	\bar{G}_{mn}(x,y,t|\tau) = \int_{0}^{L}\int_{0}^{H} G(x,y,t|\xi,\eta,\tau) \sin\mu_{m}\xi\sin\gamma_{n} \eta\dd \eta \dd \xi.
\end{equation}
Hence, $\bar{G}_{mn}$ can be viewed as an integral transform of $G$ as
in \cite{2018,2022}.

Multiplying both sides of \cref{eq:green-eq} by the function
$\sin\mu_{m}\xi\sin\gamma_{n}\eta$, and integrating over the domain
$[0,L]\times[0,H]$, we obtain a second-order differential equation for
the evolution of $\bar{G}_{mn}$, namely,
\begin{equation}
	\label{eq:green-Gmn-ode}
	\frac{\tau_{q}}{\alpha} \pdtw{\bar{G}_{mn}}{\tau} - \left(\frac{1}{\alpha} + \tau_{T} \lambda_{mn}^{2} \right) \pd{\bar{G}_{mn}}{\tau} + \lambda_{mn}^{2} \bar{G}_{mn} = \sin\mu_{m}x \sin\gamma_{n}y \,\delta(t-\tau),
\end{equation}
where $\lambda_{mn}^{2} = \mu_{m}^{2} + \gamma_{n}^{2}$. For the above
equation, it is not difficult to find that the solution, satisfying
the additional condition obtained from the integral transform of
\cref{eq:green-eq-ic}, can be expressed as
\begin{equation}
	\label{eq:green-Gmn-sol}
	\bar{G}_{mn}(x,y,t|\tau) = \frac{\alpha \sin \mu_{m}x \sin\gamma_{n} y}{\tau_{q} \beta_{2}} \exp(-\beta_{1}(t-\tau)) \sinh(\beta_{2}(t-\tau)),
\end{equation}
in which
\begin{equation}
	\label{eq:green-Gmn-sol-betas}
	\beta_{1} = \frac{1+\alpha \tau_{T} \lambda_{mn}^{2}}{2\tau_{q}}, \quad \beta_{2} = \frac{\sqrt{\left(1+\alpha\tau_{T}\lambda_{mn}^{2}\right)^{2} - 4 \alpha\tau_{q}\lambda_{mn}^{2}}}{2\tau_{q}}.
\end{equation}
Consequently, we have
\begin{equation}
	\label{eq:green-fun-series-sol}
	G(x,y,t|\xi,\eta,\tau) = \sum_{m=1}^{\infty}\sum_{n=1}^{\infty} \frac{4\alpha}{LH\tau_{q}\beta_{2}} \sin\mu_{m}x\sin\gamma_{n}y \sin\mu_{m}\xi\sin\gamma_{n}\eta \exp(-\beta_{1}(t-\tau)) \sinh(\beta_{2}(t-\tau)).  
\end{equation}

Finally, substituting
\cref{eq:green-fun-series-sol,eq:rhs-dpl-heat-eq} into
\cref{eq:integral-form-sol}, we immediately get the series expression
of the temperature as follows
\begin{equation}
	\label{eq:series-form-sol}
	T(x,y,t) = \sum_{m=1}^{\infty}\sum_{n=1}^{\infty} \frac{4\alpha\Theta}{LH\tau_{q}\beta_{2}k} P_{mn}(t) \sin\mu_{m}x\sin\gamma_{n}y,  
\end{equation}
where
\begin{align}
	\label{eq:series-form-sol-coe}
	\begin{aligned}
		P_{mn}(t) 
		= \int_{0}^{t} &\left(\sin\mu_{m} \bar{x}(\tau) \sin\gamma_{n}\bar{y}(\tau) +\tau_{q}\mu_{m} \pd{\bar{x}(\tau)}{\tau} \cos\mu_{m}\bar{x}(\tau)\sin\gamma_{n}\bar{y}(\tau)
		+ \tau_{q}\gamma_{n} \pd{\bar{y}(\tau)}{\tau} \sin\mu_{m}\bar{x}(\tau) \cos\gamma_{n}\bar{y}(\tau)\right) \\
		& \cdot \exp(-\beta_{1}(t-\tau)) \sinh(\beta_{2}(t-\tau))\dd\tau.
	\end{aligned}
\end{align}
In practice, taking the computational cost and the accuracy into
account, the infinite series of the temperature
\cref{eq:series-form-sol} would be truncated with two appropriately
selected integers $M$ and $N$ so that
\begin{equation}
	\label{eq:series-truncated-sol}
	T(x,y,t) \approx \sum_{m=1}^{M}\sum_{n=1}^{N} \frac{4\alpha\Theta}{LH\tau_{q}\beta_{2}k} P_{mn}(t) \sin\mu_{m}x\sin\gamma_{n}y.
\end{equation}
Moreover, it is obvious that the integral $P_{mn}(t)$ is too
complicated to be calculated analytically for the motion of heat
source given in \cref{eq:trajectory}. Alternatively, in the present
work, it would be calculated numerically using the \emph{Mathematica}
software \cite{Mathematica}, which is the main tool employed to
produce the results shown in the next section.
\section{Results and discussion}
\label{sec:result}

According to the previous analytical solution, temperature responses
of the plate subjected to a moving heat source with three different
trajectories, i.e., line segment trajectory (LST), circular trajectory
(CT), and elliptical trajectory (ET), are investigated in this
section. The influences of phase lags and angular velocity on the
temperature distribution are mainly analyzed. If not specified below,
the default parameters listed in \cref{tab:default-para} would be
utilized for all calculations. They are actually consistent with the
corresponding parameters employed in \cite{2008}, where temperature
distribution in a rectangular plate using the classical 3D heat
conduction model was studied.

\begin{table}[!htbp]
	\centering
	\caption{Default parameters in the calculation.}
	\label{tab:default-para}
	\begin{tabular}{cccccccccccccccc}
		\toprule
		& & & & \multicolumn{4}{c}{LST} & \multicolumn{4}{c}{CT} & \multicolumn{4}{c}{ET} \\
		\cmidrule(r){5-8}\cmidrule(r){9-12}\cmidrule(r){13-16}
		$\Theta$ & $\alpha$ & $k$ & $w$ & $L$ & $H$ & $A$ & $B$ & $L$ & $H$ & $A$ & $B$ & $L$ & $H$ & $A$ & $B$\\
		\midrule
		$2.5\times 10^{4}$ & $1.29\times 10^{-5}$ & $51.4$ & $0.2 \pi$ & $0.5$ & $0.4$ & $0.2$ & $0$ & $1.0$ & $1.0$ & $0.25$ & $0.25$ & $1.0$ & $0.5$ & $0.3$ & $0.2$ \\
		\bottomrule
	\end{tabular}
\end{table}

\subsection{Temperature response for classical heat conduction model}
\label{sec:result-classical}

As mentioned in \cref{sec:model}, the DPL heat conduction equation
\cref{eq:dpl-heat-conduction} with $\tau_{q}=\tau_{T}=0$ degenerates
to the classical heat conduction equation which is established from
Fourier's law. Apparently, however, the series
\cref{eq:series-form-sol} does not contain the solution of the
classical heat conduction model, due to the appearance of $\tau_{q}$
in the denominator. Indeed,
after a suitable modification to the solution of the second-order
equation \cref{eq:green-Gmn-ode} when it is simplified to a
first-order equation for $\tau_{q}=0$, we are able to get the solution
of this classical model as follows
\begin{equation}
	\label{eq:series-form-sol-classical}
	T(x,y,t) = \sum_{m=1}^{\infty}\sum_{n=1}^{\infty} \frac{4\alpha\Theta}{LHk} P_{mn}(t) \sin\mu_{m}x\sin\gamma_{n}y,
\end{equation}
where
\begin{align}
	\label{eq:series-form-sol-coe-classical}
	P_{mn}(t) = \int_{0}^{t} \exp(-\alpha \lambda_{mn}^{2}(t-\tau)) \sin\mu_{m} \bar{x}(\tau) \sin\gamma_{n}\bar{y}(\tau) \dd\tau.
\end{align}

\begin{figure}[!htb]
	\centering
	\includegraphics[width=0.35\textwidth]{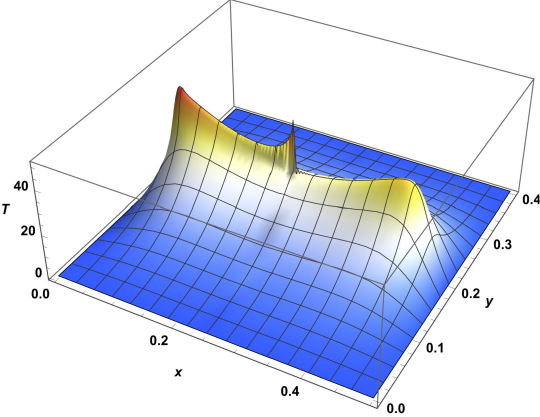}\quad
	\includegraphics[width=0.35\textwidth]{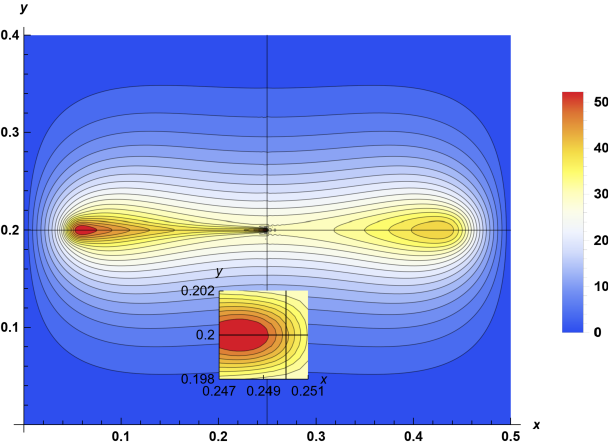}
	\caption{Temperature distribution of classical model at $t=367.5$
		for the case of LST. The intersection of horizontal and vertical
		lines denotes the location of heat source.}
	\label{fig:T-classical-lst}
\end{figure}

Although only 2D heat transfer with a moving point heat source is
under consideration in this paper, the temperature given by
\cref{eq:series-form-sol-classical} could already capture the main
features of thermal response as revealed in \cite{2008} for the cases
of LST, CT, and ET. In more detail, for the case of LST at $t=367.5$,
the temperature distribution of the classical model is shown in
\cref{fig:T-classical-lst}. At this time, the heat source is located
at the center of the plate and moving to the right with a maximum
speed of $A=0.2$. Accordingly, a temperature peak which decreases
sharply in the front and a little more gently in the rear is observed
almost at the position of the heat source. Besides, since the speed of
the heat source slows gradually to the minimum $0$ as the heat source
approaches the endpoints of the trajectory, indicating a longer
heating time in these regions, temperature peaks around both endpoints
can also be found.

\begin{figure}[!htb]
	\centering
	\subfloat[CT]{\includegraphics[width=0.3\textwidth]{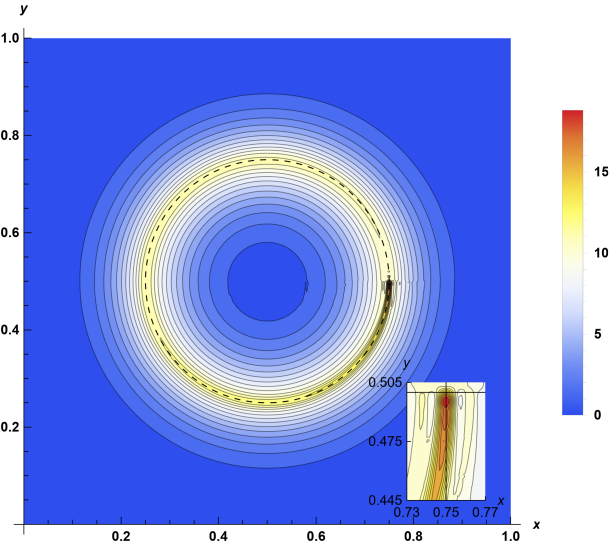}}\quad
	\subfloat[ET]{\includegraphics[width=0.54\textwidth]{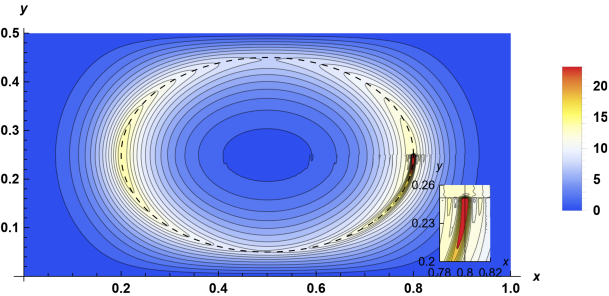}}
	\caption{Temperature distributions of classical model at $t=360$ for
		the cases of CT and ET. The intersection of horizontal and
		vertical lines denotes the location of heat source. The dashed
		curve is the trajectory of the heat source.}
	\label{fig:T-classical-ct-et}
\end{figure}

As for the cases of CT and ET, the temperature distributions of
the classical model at $t=360$ are presented in
\cref{fig:T-classical-ct-et}. As expected, the temperature along the
trajectory is much higher than other regions, and a temperature peak
with a sharp front is formed around the heat source. Meanwhile, for
the case of ET at the moment, it can be seen that, along the
trajectory, there is a temperature valley around each co-vertex
of the ellipse. This is because that the heat source moves faster at
the co-vertices than at the vertices, resulting less heat
released 
in the regions around the co-vertices.

\begin{figure}[!htb]
	\centering
	\includegraphics[width=0.6\textwidth]{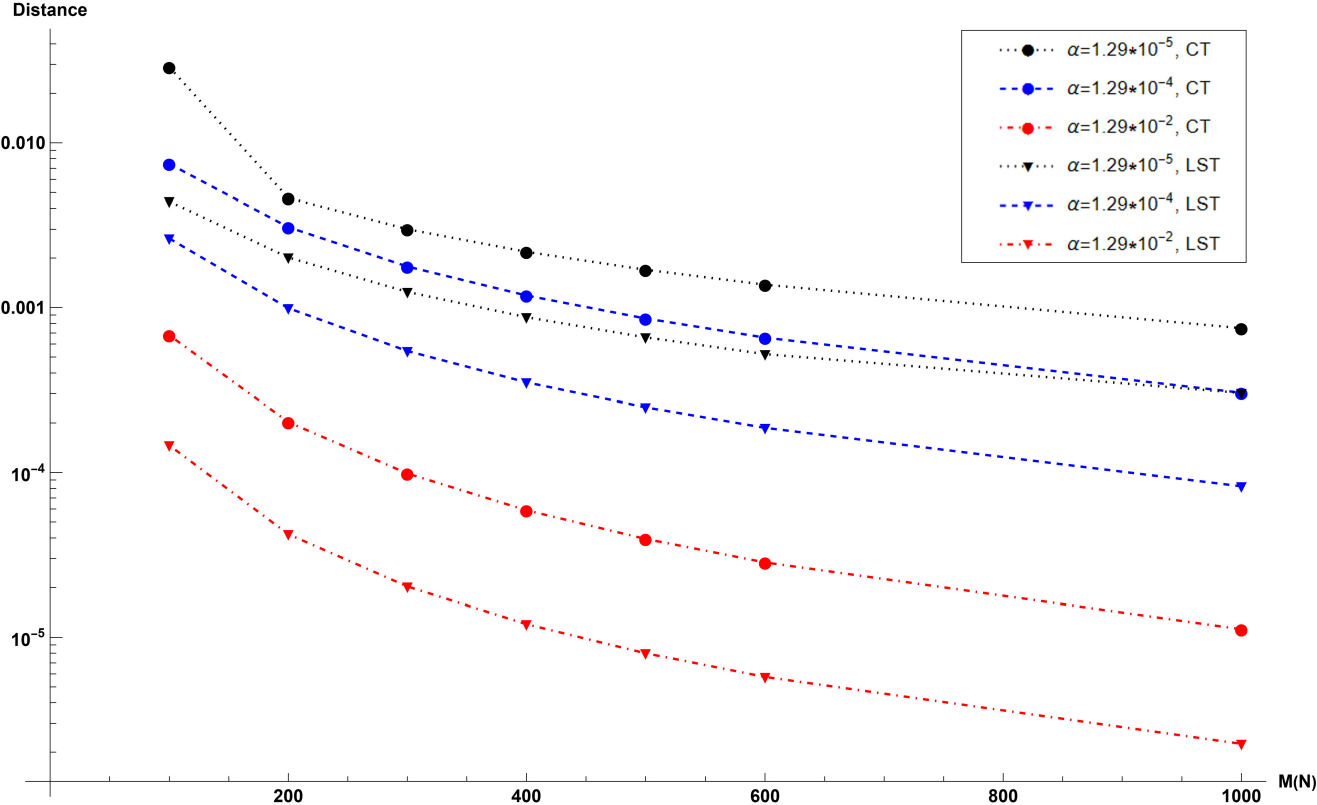}
	\caption{Distances between points of the heat source and the concomitant temperature peak in terms of the truncated integers $M(N)$.}
	\label{fig:D-classical-a}
\end{figure}

Moreover, from the zooms around the heat source in
\cref{fig:T-classical-lst,fig:T-classical-ct-et}, it is found that the
concomitant temperature peak lags slightly behind the heat source.
This phenomenon can also be observed from the results reported in
\cite{2008,hu2020heat}, where the heat source was described by other
functions rather than the Dirac delta function. It could be understood
in view of the facts that the heat source is continuously exposed to a
much cooler region, and the temperature might keep rising for a moment
after the heat source moves away, due to the rapid conduction of
considerable heat generated by the heat source.
However, in \cite{2008}, the author still stated that the highest
temperature is achieved at the heat source point.
To further investigate the position dependence of the concomitant
temperature peak on the heat source, we evaluated the distances
between them in the cases of LST and CT with three choices of
$\alpha$. The results in terms of the truncated integers $M(N)$ are
plotted in \cref{fig:D-classical-a}. It can be observed that, for all
cases, the distance decreases and would converge to $0$ as
$M,N\to \infty$. Consequently, we can deduce that for the heat source
modeled by the Dirac delta function, the temperature peak would indeed
be reached at the heat source point, although the series solution
\cref{eq:series-form-sol-classical} seems to converge a bit slowly in
the vicinity of the heat source point. In addition, we have that, with
the same truncated integers, the series solution is able to give
better result for larger $\alpha$, by noting the associated distance
is smaller.

\subsection{Temperature response for DPL heat conduction model}
\label{sec:result-dpl}

\subsubsection{DPL model for $\tau_{q}=\tau_{T}=1$ with comparison to classical  model}
\label{sec:result-dpl-equal-tau}

\begin{figure}[!htb]
	\centering
	\includegraphics[width=0.35\textwidth]{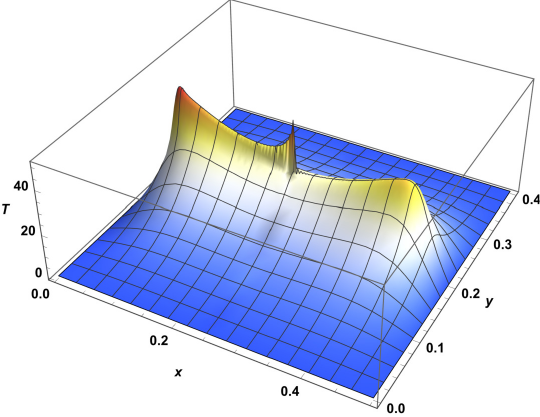}\quad
	\includegraphics[width=0.35\textwidth]{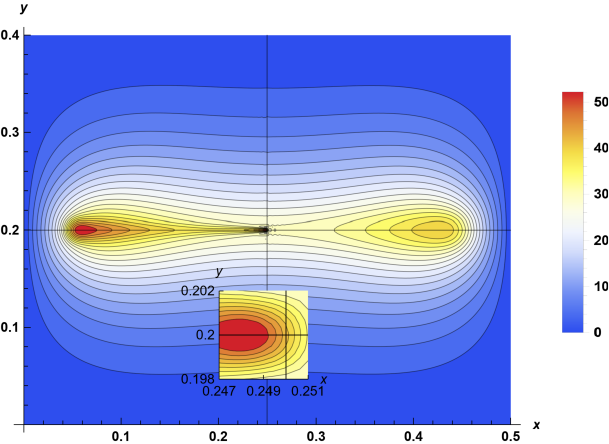}
	\caption{Temperature distribution of DPL model with
		$\tau_q=\tau_T=1$ at $t=367.5$ for the case of LST. The
		intersection of horizontal and vertical lines denotes the location
		of heat source.}
	\label{fig:T-DPL-lst-1-1}
\end{figure}

\begin{figure}[!htbp]
	\centering
	\includegraphics[width=0.6\textwidth]{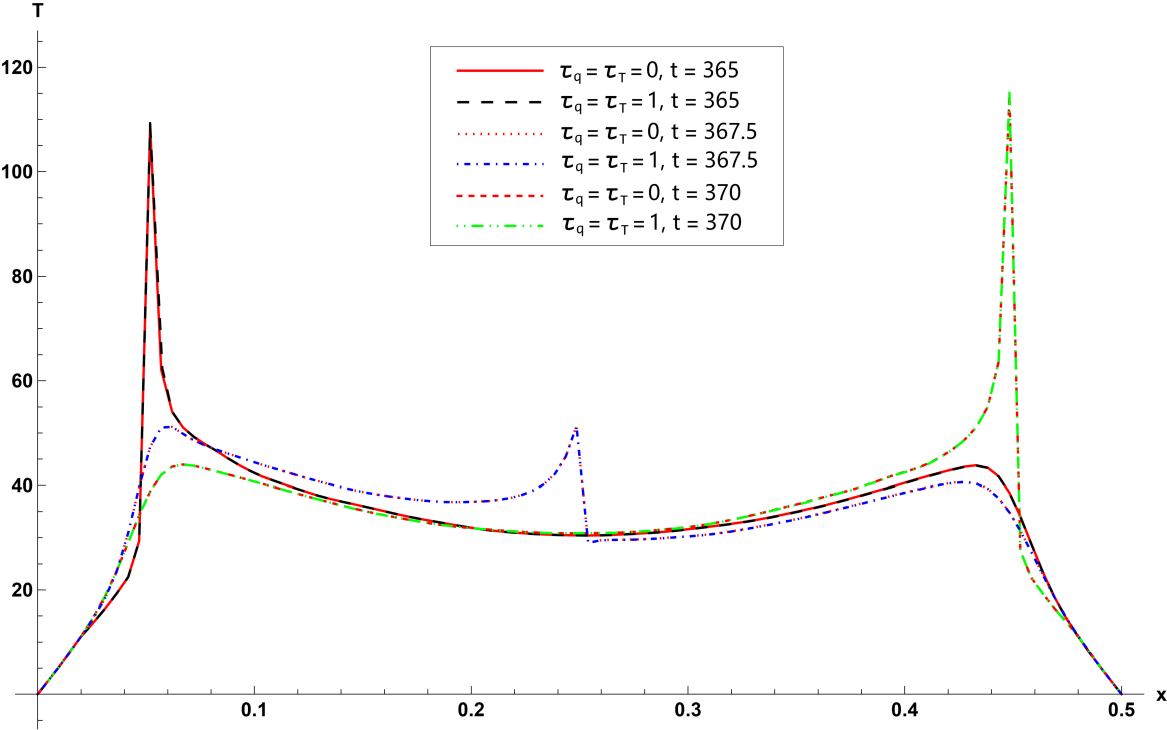}
	\caption{Temperature distributions on $y=0.2$ at three times for both
		classical model and DPL model with $\tau_{q}=\tau_{T}=1$.}
	\label{fig:T-lst-comparison-t} 
\end{figure}

\begin{figure}[!htb]
	\centering
	\subfloat[$t=362.5$]{\label{fig:T-DPL-lst-1-1-t362.5}
		\includegraphics[width=0.3\textwidth]{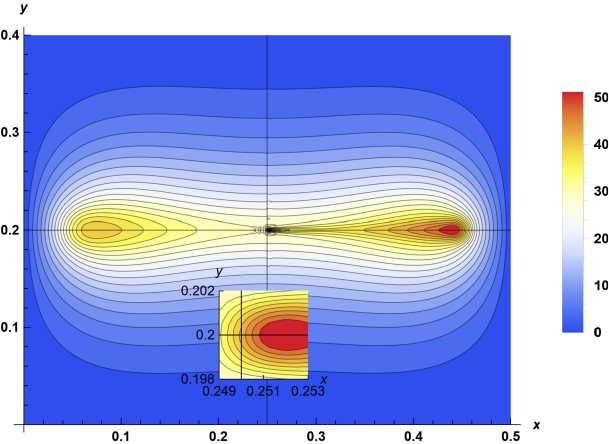}}\quad
	\subfloat[$t=365$]{\label{fig:T-DPL-lst-1-1-t365}
		\includegraphics[width=0.3\textwidth]{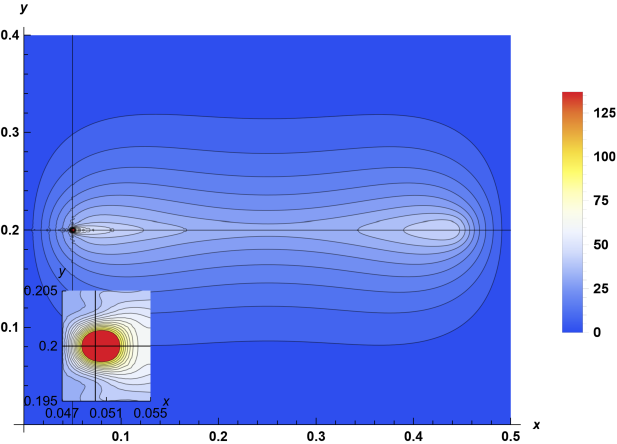}}\quad
	\subfloat[$t=370$]{\label{fig:T-DPL-lst-1-1-t370}
		\includegraphics[width=0.3\textwidth]{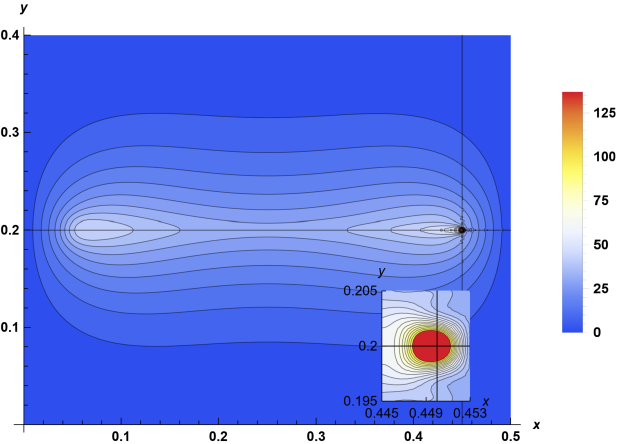}}
	\caption{Temperature distributions of DPL model with
		$\tau_q=\tau_T=1$ for the case of LST. The intersection of
		horizontal and vertical lines denotes the location of heat
		source.}
	\label{fig:T-DPL-lst-1-1-t}
\end{figure}

For the DPL heat conduction model, the phase lags are first fixed as
$\tau_{q}=\tau_{T}=1$ in the calculation. The temperature distribution
for the case of LST at $t=367.5$ is presented in
\cref{fig:T-DPL-lst-1-1}, from which nearly identical features as the
classical model are observed. Actually, in the case of
$\tau_{q}=\tau_{T}$, since the DPL constitutive relation
\cref{eq:nonFourier-law} reduces to the classical Fourier's law, it is
expected and has been validated in, e.g., \cite{2018}, that the DPL
model 
would give rise to the same results as the classical model under
uniformly distributed initial condition. This is also well verified in
\cref{fig:T-lst-comparison-t}, where, for an accurate comparison, the
temperature distributions of both models along the middle horizontal
line $y=0.2$ at three time instances are plotted. As a supplement, the
corresponding 2D temperature distributions of the DPL model are shown
in \cref{fig:T-DPL-lst-1-1-t}, which together with
\cref{fig:T-DPL-lst-1-1} exhibits the periodic variation
characteristic of the temperature distribution. Particularly, the
temperature distribution at $t=365$ is symmetrical with that at
$t=370$. So are the temperature distributions at $t=362.5$ and
$t=367.5$.

\begin{figure}[!htb]
	\centering
	\subfloat[CT]{\label{fig:T-DPL-ct-1-1}
		\includegraphics[width=0.3\textwidth]{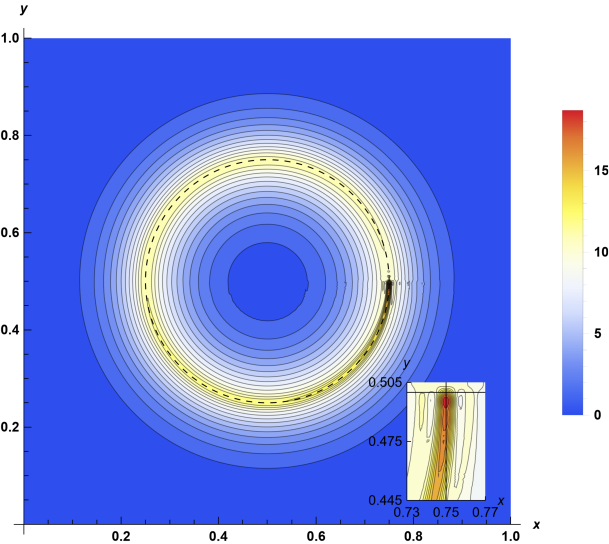}}\quad
	\subfloat[ET]{\label{fig:T-DPL-et-1-1}
		\includegraphics[width=0.54\textwidth]{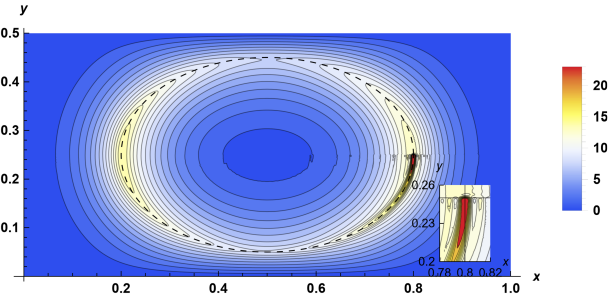}}
	\caption{Temperature distributions of DPL model with
		$\tau_{q}=\tau_{T}=1$ at $t=360$ for the cases of CT and ET. The
		intersection of horizontal and vertical lines denotes the location
		of heat source. The dashed curve is the trajectory of heat
		source.}
	\label{fig:T-DPL-ct-et-1-1}
\end{figure}

Similarly, for the DPL model with $\tau_{q}=\tau_{T}=1$ in the cases
of CT and ET, it is shown in \cref{fig:T-DPL-ct-et-1-1} as well as
\cref{fig:T-classical-ct-et} that the same results as the classical
model would be obtained too.



\subsubsection{Influences of phase lags}
\label{sec:result-dpl-neq-tau}

As illustrated in many studies, e.g., \cite{tzou2014macro, 2018,2022},
both phase lags $\tau_{q}$ and $\tau_{T}$ play important roles in the
DPL model. When $\tau_{q}\neq \tau_{T}$, the resulting temperature
distribution exhibits a quite different characteristic from the
classical model. Meanwhile, for the heat source moving along LST, CT,
and ET, the DPL model would give rise to a distinctive temperature
response, in comparison to the common situation that the heat source
moves along a straight line with a constant velocity.

\begin{figure}[!htb]
	\centering
	\subfloat[$\tau_{q}=1,\tau_{T}=5$]{\label{fig:T-DPL-lst-1-5}
		\includegraphics[width=0.35\textwidth]{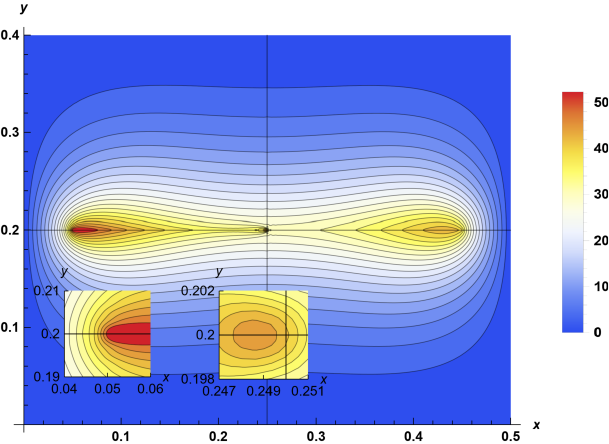}}\quad
	\subfloat[$\tau_{q}=5,\tau_{T}=1$]{\label{fig:T-DPL-lst-5-1}
		\includegraphics[width=0.35\textwidth]{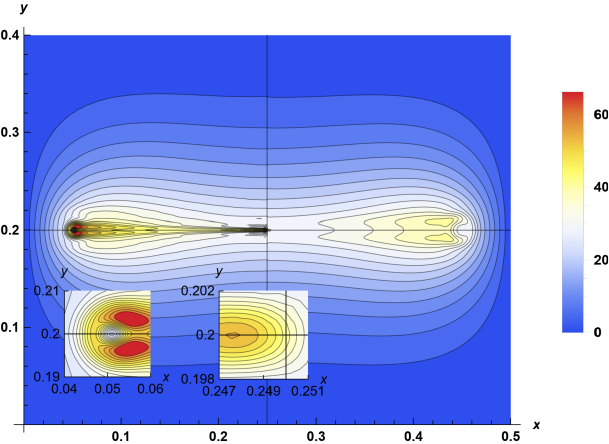}}
	\caption{Temperature distributions of DPL model with two pairs of
		phase lags at $t=367.5$ for the case of LST. The intersection of horizontal and
		vertical lines denotes the location of heat source.}
	\label{fig:T-DPL-lst}
\end{figure}

\begin{figure}[!htb]
	\centering
	\subfloat[$t=365$]{\label{fig:T-DPL-lst-5-1-t365}
		\includegraphics[width=0.3\textwidth]{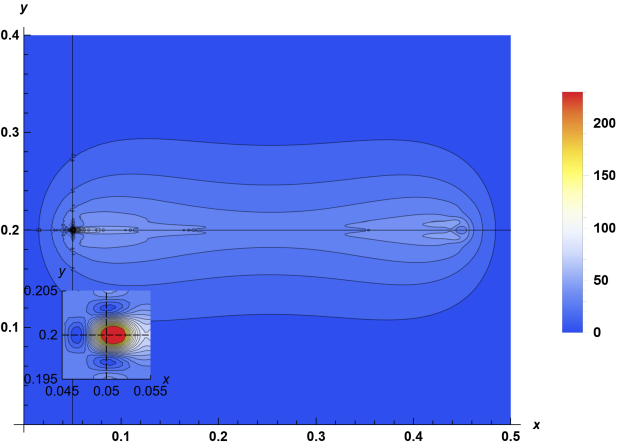}}\quad
	\subfloat[$t=365.5$]{\label{fig:T-DPL-lst-5-1-t365.5}
		\includegraphics[width=0.3\textwidth]{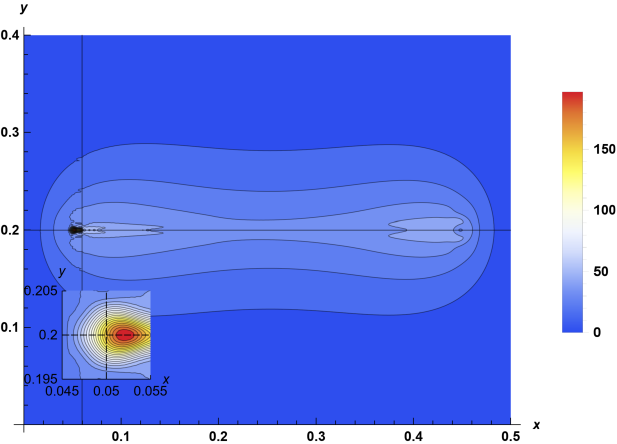}}\quad
	\subfloat[$t=366$]{\label{fig:T-DPL-lst-5-1-t366}
		\includegraphics[width=0.3\textwidth]{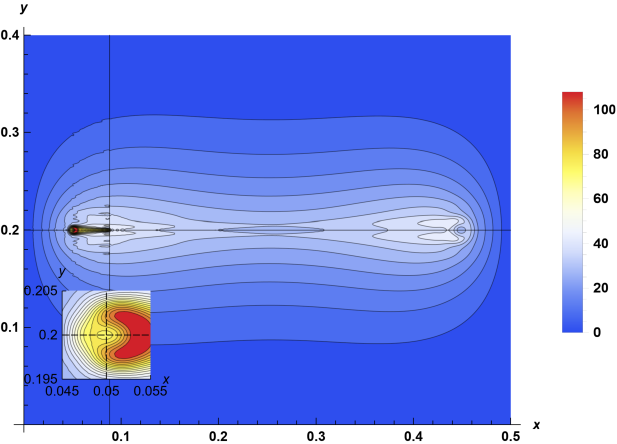}}
	\caption{Temperature distributions of DPL model with $\tau_q=5$,
		$\tau_T=1$ for the case of LST. The intersections of solid lines
		and dashed lines in the horizontal and vertical directions
		represent, respectively, the location of heat source and the left
		endpoint of trajectory.}
	\label{fig:T-DPL-lst-5-1-t}
\end{figure}

In \cref{fig:T-DPL-lst}, the temperature distributions of the DPL
model with two pairs of phase lags given by $\tau_{q}=1$,
$\tau_{T}=5$, and $\tau_{q}=5$, $\tau_{T}=1$ at $t=367.5$ for the case
of LST are presented. Compared with the results of the classical model
shown in \cref{fig:T-classical-lst}, a lower temperature peak around
the heat source is found in \cref{fig:T-DPL-lst-1-5}, which is
reasonable as the heat flux precedes the temperature gradient such
that the temperature would be homogenized more efficiently when
$\tau_{q}<\tau_{T}$. Moreover, due to the comprehensive effect of the
tardy-temperature-rise phenomenon as observed in \cite{2018} and a
longer irradiation time of heat source at both endpoints of the
trajectory than the center, a larger and higher temperature zone
appears near the left endpoint of the trajectory. In contrast, when
$\tau_{q}>\tau_{T}$, as shown in \cref{fig:T-DPL-lst-5-1}, the
temperature gradient precedes the heat flux, and the DPL model would
embody a distinct fluctuation feature. Consequently, a higher maximum
temperature is obtained and a shallow furrow, which splits the
temperature distribution into two symmetrical parts, is detected
exactly along the trajectory of heat source. Apart from this, it is
astonished that a deep temperature well can be observed at the left
endpoint of the trajectory. To investigate it in more detail, we plot
the temperature distributions at three earlier time instants in
\cref{fig:T-DPL-lst-5-1-t}. It can be seen that at $t=365$ when the
heat source is located at the left endpoint of the trajectory, a
single high temperature zone arises nearby the heat source with the
peak much higher than that of the classical model shown in
\cref{fig:T-DPL-lst-1-1-t365}. As the heat source moves away from the
endpoint to the center with increasing speed, the concomitant
temperature peak, which closely follows the heat source, is separated
from the high temperature zone, leading to a swift decrease of the
maximum temperature. What is more, because of the volatility in heat
conduction, the high temperature zone nearby the left endpoint is
gradually split into two symmetrical high temperature zones, and the
temperature well at the left endpoint is simultaneously evolved.


\begin{figure}[!htb]
	\centering
	\subfloat[CT]{\label{fig:T-DPL-ct-1-5}
		\includegraphics[width=0.3\textwidth]{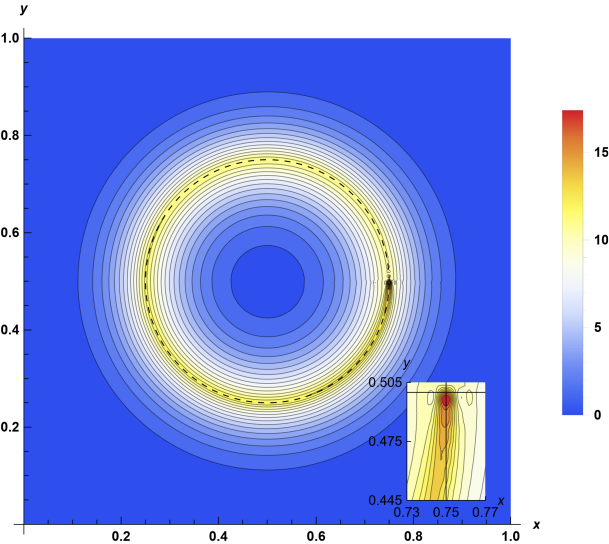}}\quad
	\subfloat[ET]{\label{fig:T-DPL-et-1-5}
		\includegraphics[width=0.54\textwidth]{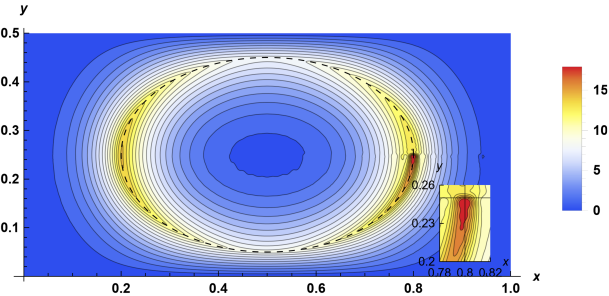}}
	\caption{Temperature distributions of DPL model with $\tau_{q}=1$,
		$\tau_{T}=5$ at $t=360$ for the cases of CT and ET. The
		intersection of horizontal and vertical lines denotes the location
		of heat source. The dashed curve is the trajectory of heat
		source.}
	\label{fig:T-DPL-ct-et-1-5}
\end{figure}

\begin{figure}[!htb]
	\centering
	\subfloat[CT]{\label{fig:T-DPL-ct-5-1}
		\includegraphics[width=0.3\textwidth]{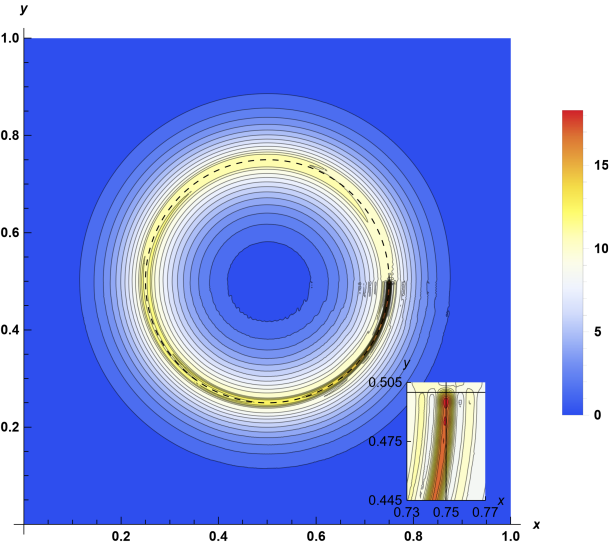}}\quad
	\subfloat[ET]{\label{fig:T-DPL-et-5-1}
		\includegraphics[width=0.54\textwidth]{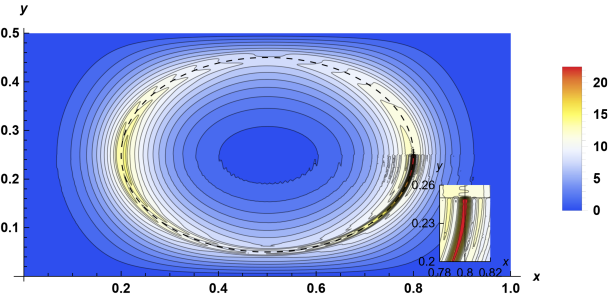}}
	\caption{Temperature distributions of DPL model with $\tau_{q}=5$,
		$\tau_{T}=1$ at $t=360$ for the cases of CT and ET. The
		intersection of horizontal and vertical lines denotes the location
		of heat source. The dashed curve is the trajectory of heat
		source.}
	\label{fig:T-DPL-ct-et-5-1}
\end{figure}

For the cases of CT and ET, the temperature distributions of the DPL
model at $t=360$ are shown in
\cref{fig:T-DPL-ct-et-1-5,fig:T-DPL-ct-et-5-1}, respectively, for
$\tau_{q}=1$, $\tau_{T}=5$ and $\tau_{q}=5$, $\tau_{T}=1$. As
mentioned previously, the precedence of heat flux for
$\tau_{q}<\tau_{T}$ implies a more efficient heat conduction process
than that for $\tau_{q}=\tau_{T}$. It turns out that a lower
temperature peak following the heat source and a gentler temperature
variation along the trajectory than the results shown in
\cref{fig:T-DPL-ct-et-1-1} are observed in
\cref{fig:T-DPL-ct-et-1-5}. Similarly, owing to the fluctuation
characteristic for $\tau_{q}>\tau_{T}$, an inconspicuous furrow along
each trajectory would also appear in \cref{fig:T-DPL-ct-et-5-1}.

\begin{figure}[!htb]
	\centering
	\subfloat[classical $\tau_{q}=\tau_{T}=0$]{\label{fig:T-DPL-a2-ct_0-0}
		\includegraphics[width=0.24\textwidth]{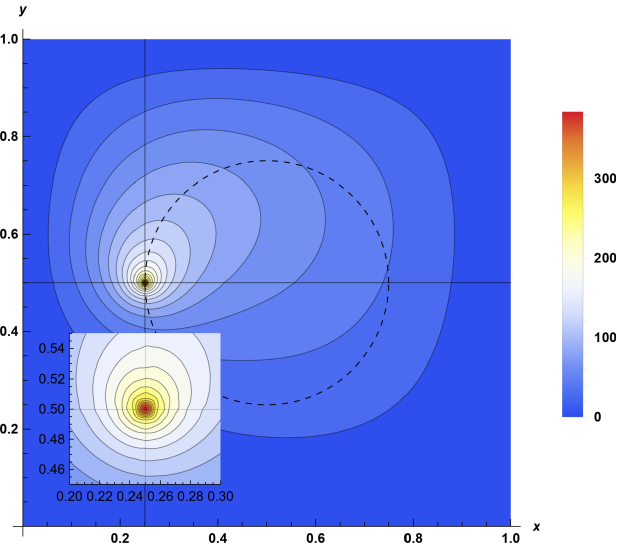}}~
	\subfloat[$\tau_{q}=1, \tau_{T}=2$]{\label{fig:T-DPL-a2-ct_1-2}
		\includegraphics[width=0.24\textwidth]{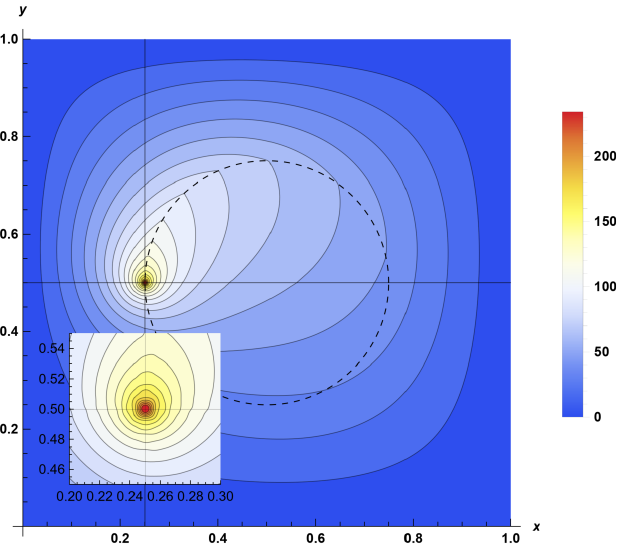}}~
	\subfloat[$\tau_{q}=1, \tau_{T}=5$]{\label{fig:T-DPL-a2-ct_1-5}
		\includegraphics[width=0.24\textwidth]{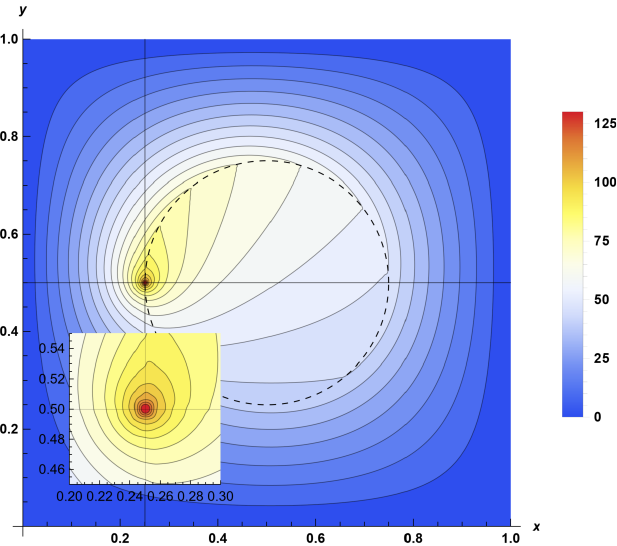}}~
	\subfloat[$\tau_{q}=1, \tau_{T}=10$]{\label{fig:T-DPL-a2-ct_1-10}
		\includegraphics[width=0.24\textwidth]{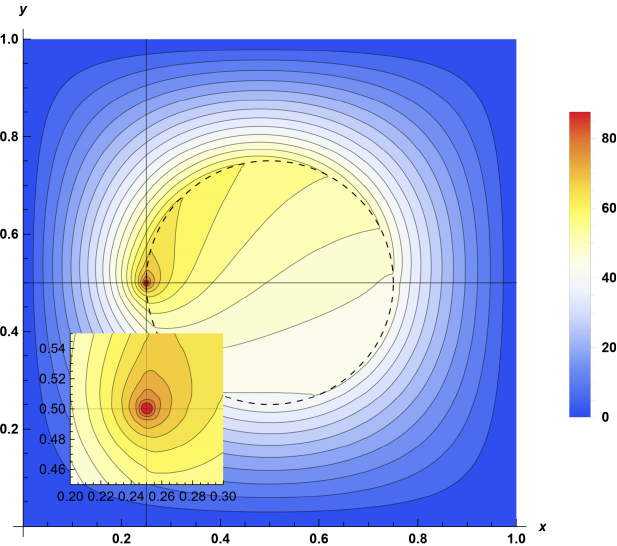}}\\
	\subfloat[$\tau_{q}=\tau_{T}=1$]{\label{fig:T-DPL-a2-ct_1-1}
		\includegraphics[width=0.24\textwidth]{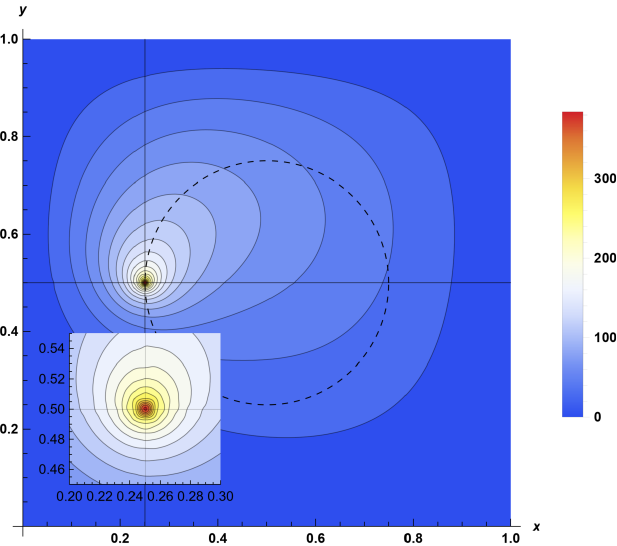}}~
	\subfloat[$\tau_{q}=2, \tau_{T}=1$]{\label{fig:T-DPL-a2-ct_2-1}
		\includegraphics[width=0.24\textwidth]{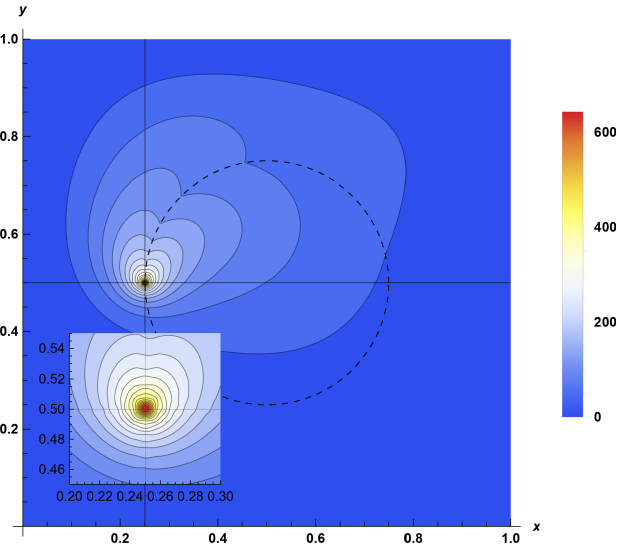}}~
	\subfloat[$\tau_{q}=5, \tau_{T}=1$]{\label{fig:T-DPL-a2-ct_5-1}
		\includegraphics[width=0.24\textwidth]{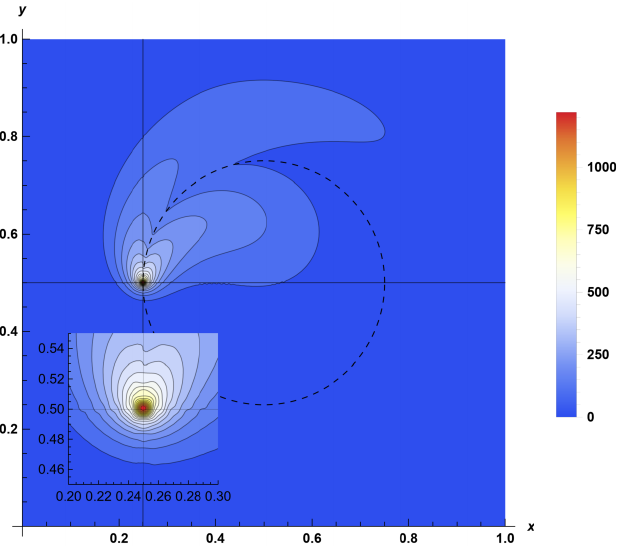}}~
	\subfloat[$\tau_{q}=10, \tau_{T}=1$]{\label{fig:T-DPL-a2-ct_10-1}
		\includegraphics[width=0.24\textwidth]{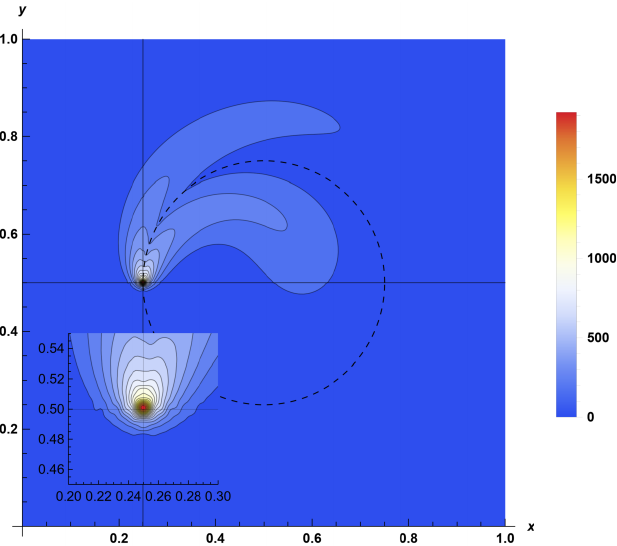}}
	\caption{Temperature distributions of DPL model with
		$\alpha=1.29\times 10^{-2}$ and various phase lags at $t=25$ for
		the case of CT. The intersection of horizontal and vertical lines
		denotes the location of heat source. The dashed curve is the
		trajectory of heat source.}
	\label{fig:T-DPL-a2-ct}
\end{figure}

It is noteworthy that, for all the above experiments, the temperature
rise is concentrated in a narrow band centered on the trajectory due
to the small thermal diffusivity $\alpha$. In order to more
expediently explore the detailed effects of both phase lags $\tau_{q}$
and $\tau_{T}$ on the temperature response, below we adjust the value
of $\alpha$ to $1.29\times 10^{-2}$, which not only enhances the
thermal diffusion capability, but also in some sense increases the
strength of the heat source. For simplicity, only the results for the
case of CT are presented in the following. In particular,
\cref{fig:T-DPL-a2-ct} shows the temperature distributions of the DPL
model with various pairs of phase lags at $t=25$.
Obviously, it can be confirmed again from the first column of
\cref{fig:T-DPL-a2-ct} that taking the same value for both phase lags,
i.e., $\tau_{q}=\tau_{T}=1$, the DPL model would give rise to the
consistent temperature distribution as the classical model with $\tau_{q}=\tau_{T}=0$.
It can also be seen from \cref{fig:T-DPL-a2-ct} that for all pairs of
phase lags a steep temperature peak with a sharp front and a
relatively gentle rear appears in the region, where the heat mainly
accumulates due to the direct heating of the heat source. In more
detail, when $\tau_{q}<\tau_{T}$ and $\tau_{q}=1$, it is anticipated
that the heat flux precedes the temperature gradient more, indicating
a more efficient diffusion, as the phase lag $\tau_{T}$
increases. Hence, from the first row of \cref{fig:T-DPL-a2-ct}, we
have that the temperature peak decreases and the temperature variation
becomes more gradual for a larger $\tau_{T}$. Besides, owing to the
limited inner region of circular trajectory,
therein the heat would accumulate more easily than the outer region of
the trajectory. Accordingly, the average temperature rise is higher
and the temperature variation is more gradual in the inner region than
that in the outer region. Moreover, it is interesting to see that the
temperature distribution in the inner and outer regions are separated
exactly by the trajectory, where the temperature ridgeline becomes
more and more evident as $\tau_{T}$ increases. On the other hand,
when $\tau_{q}>\tau_{T}$, the volatility appears in the DPL model,
resulting in that the heat shall transfer in a wave form with a finite
velocity. As mentioned in \cite{2018}, the wave velocity would be
given by $c=\sqrt{\alpha/\tau_{q}}$ for $\tau_{T}=0$. So it can be
expected to decrease as $\tau_{q}$ increases for $\tau_{T}\neq
0$. Subsequently, as presented in the second row of
\cref{fig:T-DPL-a2-ct} where $\tau_{T}=1$, we have that more heat
accumulates in the directly heated position as $\tau_{q}$ increases,
leading to a much higher temperature peak around the heat source with
a sharper front for a larger value of $\tau_{q}$. Concurrently, the
furrow, that manifests exactly along the trajectory and splits the
temperature distribution into the inner and outer branches, becomes
more conspicuous and deeper as $\tau_{q}$ increases. In addition,
slenderer contours for both branches with slightly fatter inner
contour, resembling long curved ears of a rabbit, are found when
$\tau_{q}$ enlarges.

\begin{figure}[!htb]
	\centering
	\subfloat[$\tau_{q}=1$]{\label{fig:T-DPL-a2-ct-comparison_q}
		\includegraphics[width=0.45\textwidth]{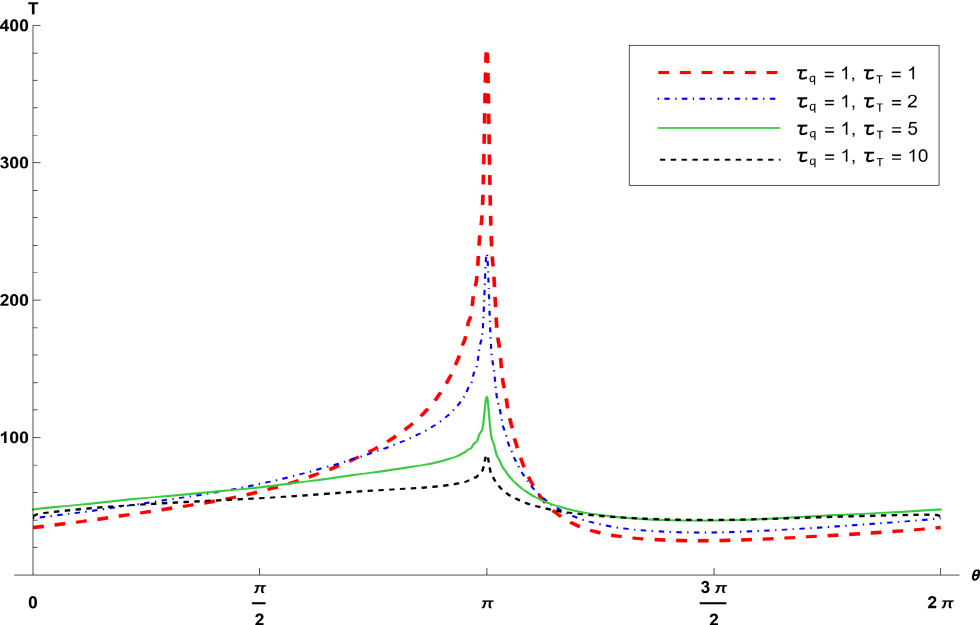}}\quad
	\subfloat[$\tau_{T}=1$]{\label{fig:T-DPL-a2-ct-comparison_T}
		\includegraphics[width=0.45\textwidth]{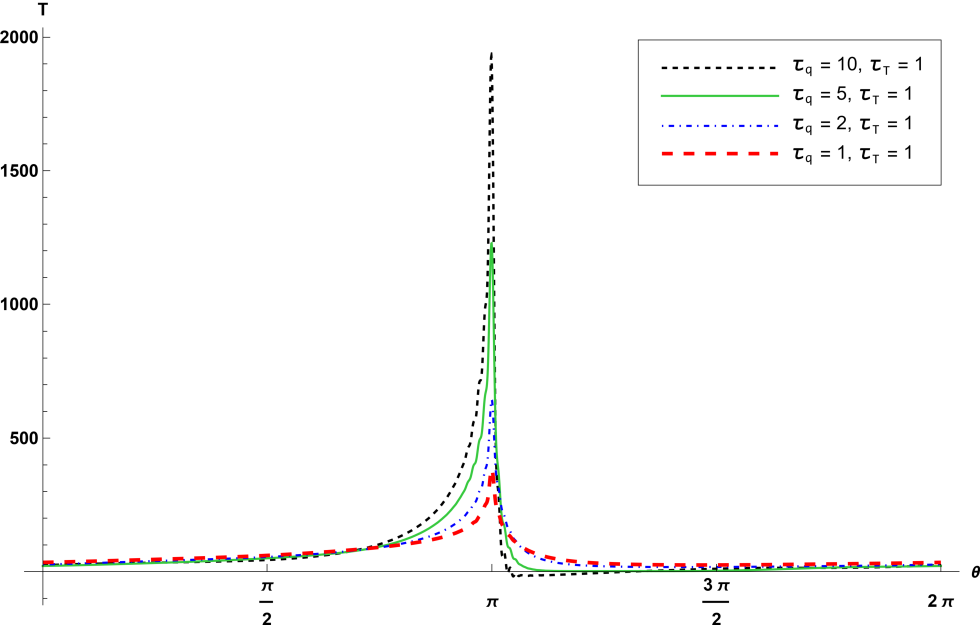}}
	\caption{Temperature distributions along CT for the DPL model with
		$\alpha=1.29\times 10^{-2}$ and various phase lags at $t=25$.}
	\label{fig:T-DPL-a2-ct-1D}
\end{figure}

To delve into the temperature response of the DPL model along the
trajectory, we plot the one-dimensional temperature distributions
along the trajectory with respect to the counterclockwise central
angle at $t=25$ for various pairs of phase lags in
\cref{fig:T-DPL-a2-ct-1D}, where the results illustrate the previous
qualitative feature again. That is, a larger $\tau_{T}$ shall give
rise to a lower temperature peak and a smaller temperature variation,
while a larger $\tau_{q}$ would lead to a much higher temperature
peak. Furthermore, as $\tau_{q}$ increases with fixed $\tau_{T}$, the
wave velocity of heat reduces, indicating that the heat generated by
the direct heating of the heat source is more difficult to transfer to
a distant position. Therefore, a steeper front of the temperature peak
appears when $\tau_{q}$ grows, and the temperature away from the heat
source is affected about the same for different $\tau_{q}$ when it is
larger than $\tau_{T}=1$.

\subsubsection{Influences of angular velocity}
\label{sec:result-dpl-var-vel}

\begin{figure}[!htb]
	\centering
	\subfloat[$w=0.1\pi$]{\label{fig:T-DPL-a2-ct-1-5-w0.1pi}
		\includegraphics[width=0.3\textwidth]{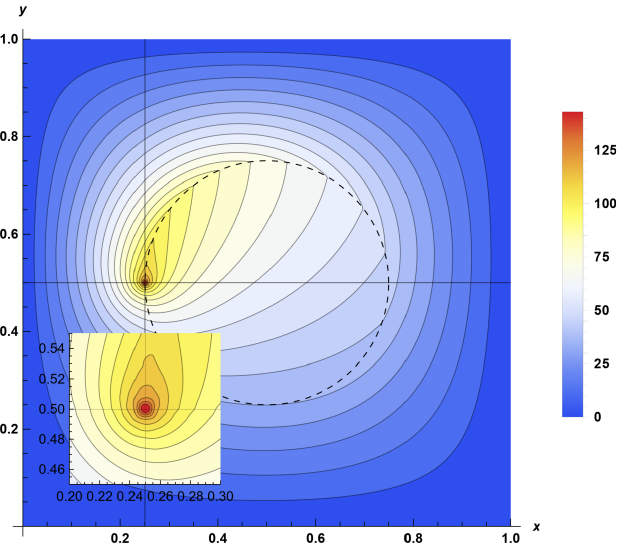}}~
	\subfloat[$w=0.2\pi$]{\label{fig:T-DPL-a2-ct-1-5-w0.2pi}
		\includegraphics[width=0.3\textwidth]{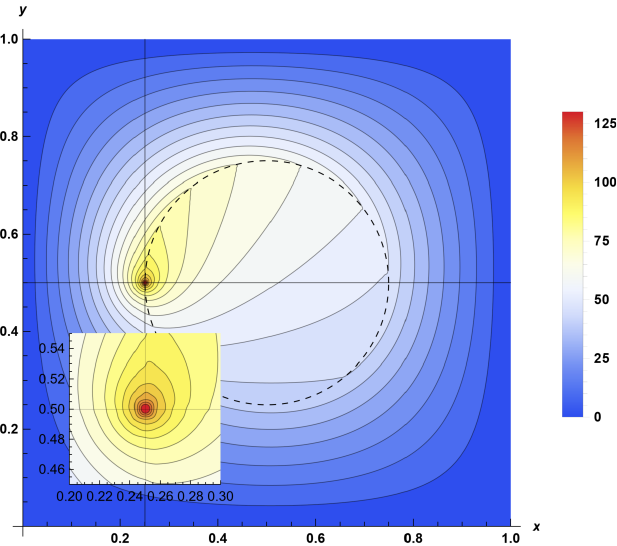}}~
	\subfloat[$w=0.4\pi$]{\label{fig:T-DPL-a2-ct-1-5-w0.4pi}
		\includegraphics[width=0.3\textwidth]{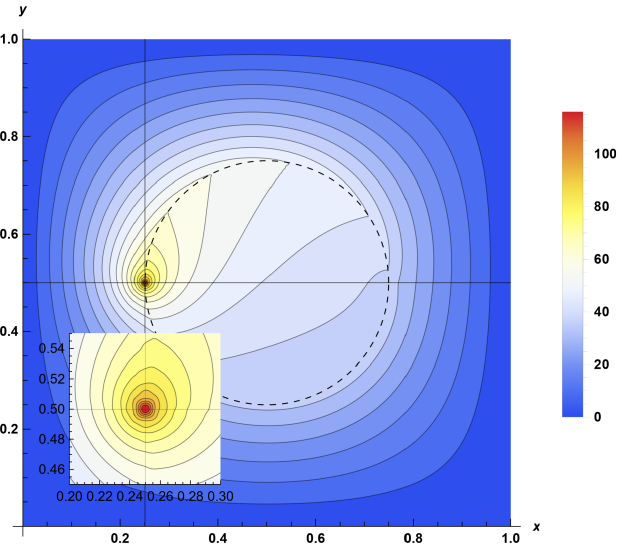}}\\
	\subfloat[$w=0.1\pi$]{\label{fig:T-DPL-a2-ct-5-1-w0.1pi}
		\includegraphics[width=0.3\textwidth]{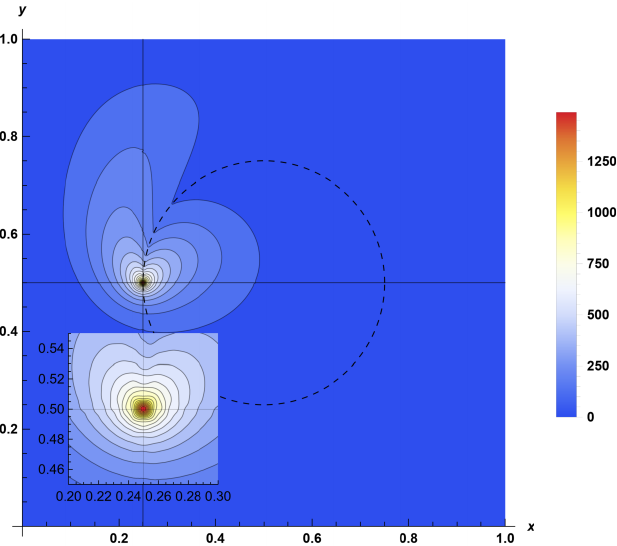}}~
	\subfloat[$w=0.2\pi$]{\label{fig:T-DPL-a2-ct-5-1-w0.2pi}
		\includegraphics[width=0.3\textwidth]{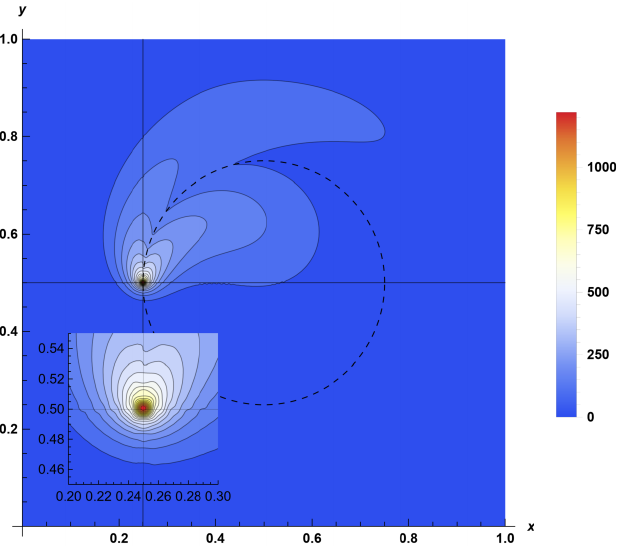}}~
	\subfloat[$w=0.4\pi$]{\label{fig:T-DPL-a2-ct-5-1-w0.4pi}
		\includegraphics[width=0.3\textwidth]{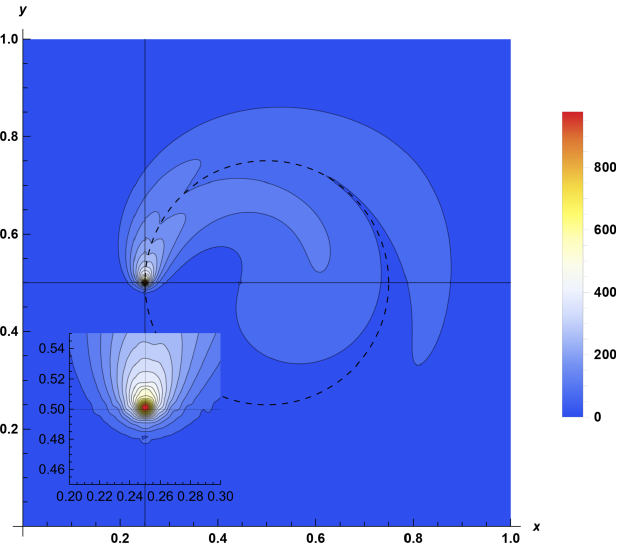}}
	\caption{Temperature distributions of DPL model with
		$\alpha=1.29\times 10^{-2}$ and various angular velocities at
		$t=70$ for the case of CT, when $\tau_{q}=1$, $\tau_{T}=5$ (top
		row) and $\tau_{q}=5$, $\tau_{T}=1$ (bottom row). The intersection
		of horizontal and vertical lines denotes the location of heat
		source. The dashed curve is the trajectory of heat source.}
	\label{fig:T-DPL-a2-ct-w}
\end{figure}

Apparently, the angular velocity $w$ is another important factor
affecting the temperature response. Here we fix the pair of phase lags
as $\tau_{q}=1$, $\tau_{T}=5$, and $\tau_{q}=5$, $\tau_{T}=1$,
respectively. The resulting temperature distributions of the DPL model
with three different angular velocities, namely, $w=0.1\pi$, $0.2\pi$,
and $0.4\pi$, are displayed in \cref{fig:T-DPL-a2-ct-w}. In view of
that a smaller angular velocity implies a longer irradiation time for
a given region along the trajectory, we can expect and observe in
\cref{fig:T-DPL-a2-ct-w} that the concomitant temperature peak, which
is mainly caused by direct heating, is higher when the angular
velocity is smaller for both pairs of phase lags. Moreover, a faster
angular velocity means a shorter motion period, which also connotes
that a given region along the trajectory would be directly heated by
the heat source more frequently. It follows that, under the combined
effects of thermal diffusion and moving heating in the situation the
heat flux precedes the temperature gradient, i.e.,
$\tau_{q}<\tau_{T}$, the temperature away from the heat source,
especially in the inner region of the trajectory, varies more gently
when the angular velocity is larger. Conversely, in the situation of
$\tau_{q}>\tau_{T}$, it can be seen from the bottom row of
\cref{fig:T-DPL-a2-ct-w} that the temperature contours, like long
curved ears of a rabbit, become slenderer as the angular velocity
increases.

\begin{figure}[!htb]
	\centering
	\subfloat[$\tau_{q}=1, \tau_{T}=5$]{\label{fig:T-DPL-a2-ct-comparison_1-5}
		\includegraphics[width=0.45\textwidth]{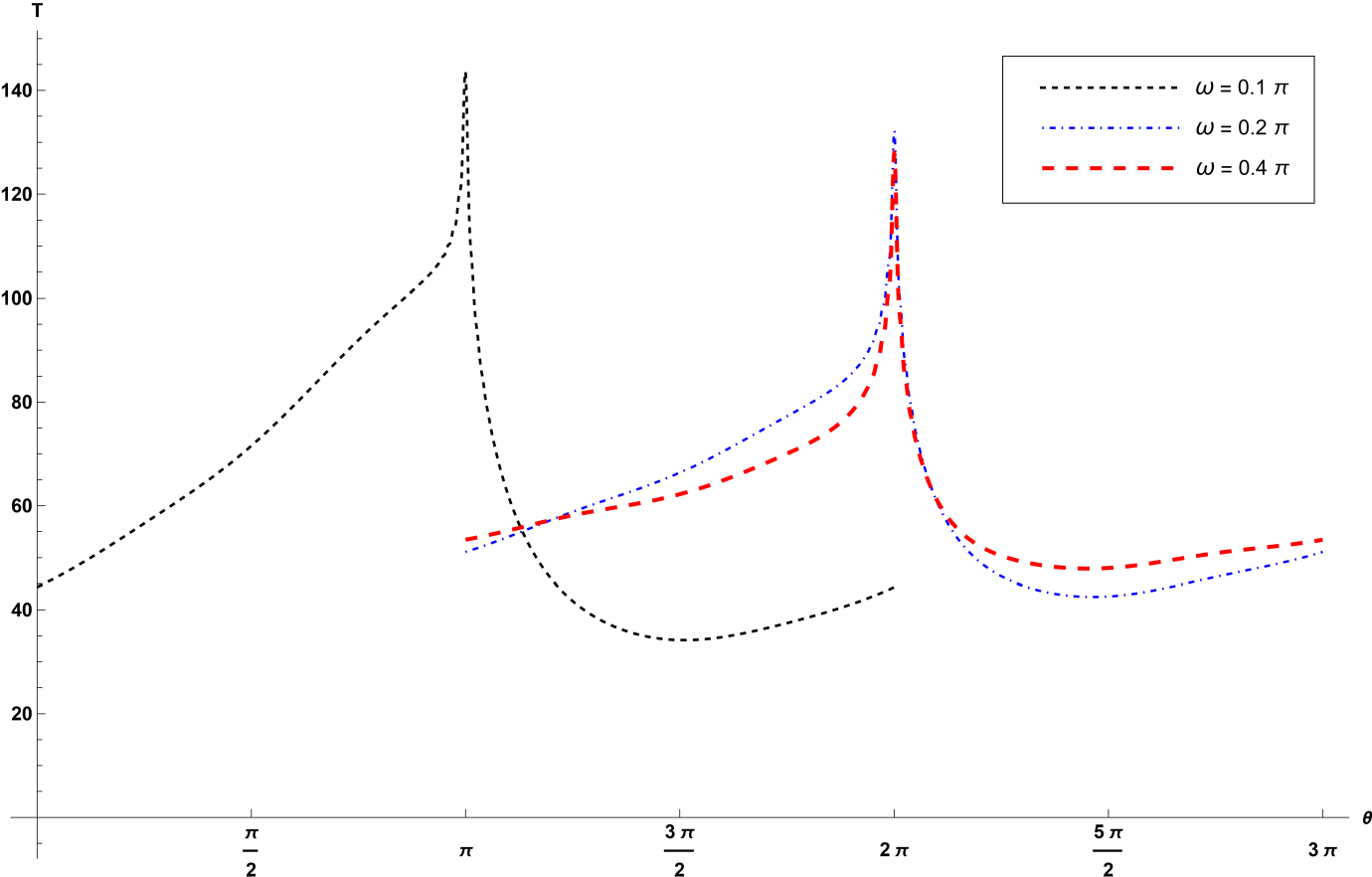}}\quad
	\subfloat[$\tau_{q}=5, \tau_{T}=1$]{\label{fig:T-DPL-a2-ct-comparison_5-1}
		\includegraphics[width=0.45\textwidth]{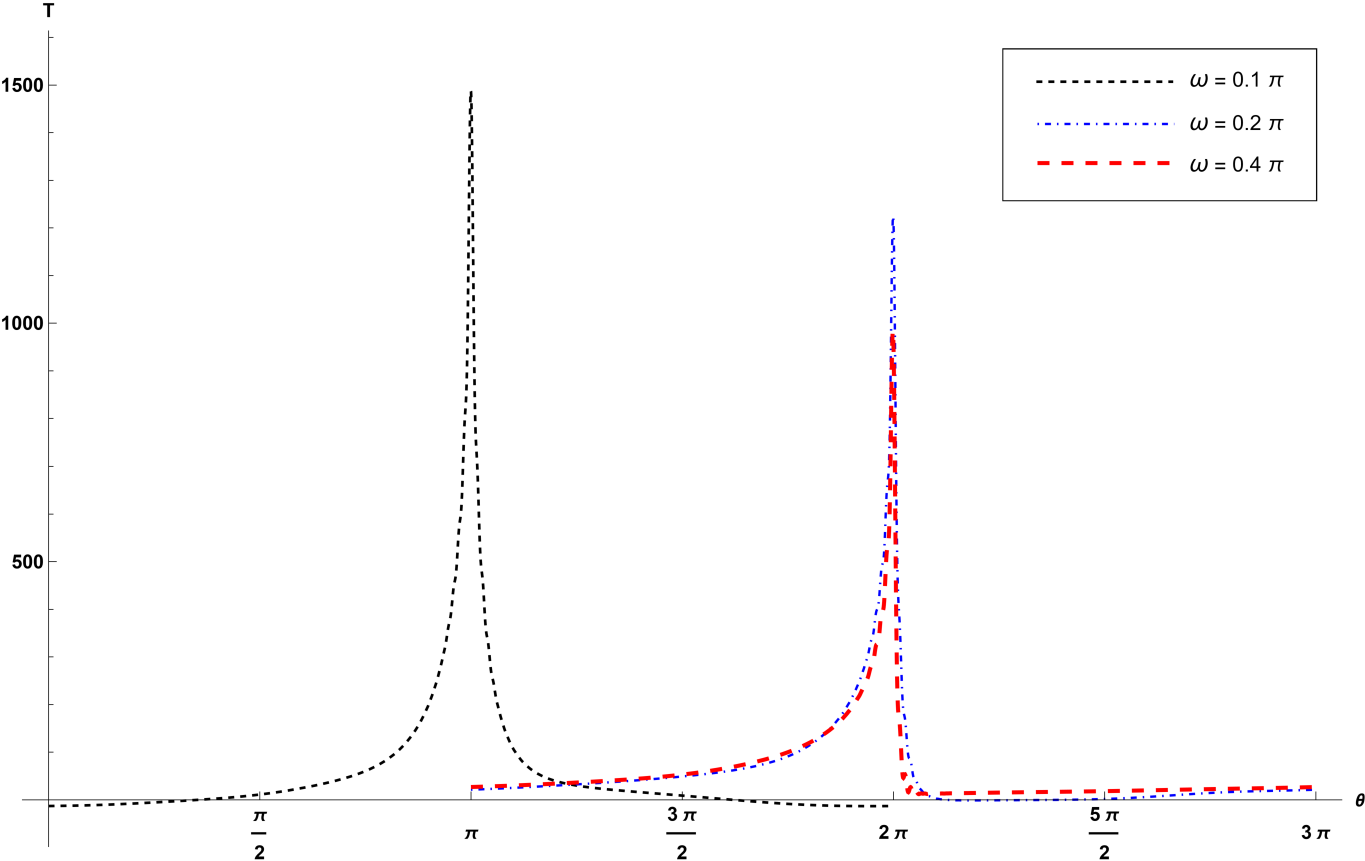}}
	\caption{Temperature distributions along CT for the DPL model with
		$\alpha=1.29\times 10^{-2}$ and various angular velocities at $t=70$.}
	\label{fig:T-DPL-a2-ct-1D-w}
\end{figure}

At last, the corresponding one-dimensional temperature distributions
along the trajectory with respect to the counterclockwise central
angle are shown in \cref{fig:T-DPL-a2-ct-1D-w}. Again, it can be
observed that the temperature peak is lower and the temperature away
from the heat source varies more gently, when the angular velocity is
larger for both pairs of phase lags. Apart from this, in the situation
of $\tau_{q}<\tau_{T}$, it can also be seen that the height of the
temperature peak for $w=0.2\pi$ is very close to that for
$w=0.4\pi$. This might be understood that a delicate balance has been
reached between thermal diffusion and moving heating. Whereas in the
situation of $\tau_{q}>\tau_{T}$, since the heat is difficult to
transfer to a distant position, we can find in
\cref{fig:T-DPL-a2-ct-comparison_5-1} that a larger angular velocity
would lead to a sharper front of the temperature peak.

\section{Conclusion}
\label{sec:conclusion}

The DPL heat conduction process in a two-dimensional plate, irradiated
by a point heat source that moves along three different trajectories,
namely, LST, CT, and ET, has been studied in this work. The analytical
solution of temperature for the underlying heat conduction model is
first derived by employing Green's function approach. Based on the
series representation of this analytical solution, the thermal
responses for the present moving heat source problem are then
investigated. The relations between the moving heat source and the
concomitant temperature peak, the influences of the pair of phase lags
and the angular velocity on the temperature are mainly analyzed and
discussed. The results show quite distinctive and different thermal
behaviors for all three trajectories, in
comparison to the results revealed for the common situation that the
heat source moves in a straight line with a constant speed. To be
specific, in the DPL model, the appearance of the phase lag $\tau_{q}$
brings in fluctuation feature to the heat transfer, and the
introduction of the phase lag $\tau_{T}$ makes the heat tend to spread
more efficiently. When $\tau_{q}<\tau_{T}$, the heat diffusion
dominates, and its efficiency is improved more and more as $\tau_{T}$
increases, resulting in a larger region of temperature rise and a
lower temperature peak. When $\tau_{q}>\tau_{T}$, the fluctuation of
heat transfer dominates, and is more and more evident, accompanied
with a decreasing velocity of heat wave, as $\tau_{q}$ increases. It
follows that, for a larger $\tau_{q}$, the heat is more difficult to
transfer to a long distance, leading to a much higher temperature peak
around the heat source. Along the trajectory of the heat source, a
shallow furrow that splits the temperature distribution into two parts
can also be observed. In the case of LST, a deep temperature well is
even evolved at the endpoint of the trajectory after the heat source
moves away from it. For the remaining situation $\tau_{q}=\tau_{T}$,
the DPL model gives rise to the same results as the
classical model.
Furthermore, with a large thermal diffusivity, the heat can spread
more easily. In the case of CT, we subsequently observe that the
temperature distribution in the inner and outer regions of the
trajectory are separated clearly for all situations except for
$\tau_{q}=\tau_{T}$. Precisely speaking, there exists a ridgeline for
$\tau_{q}<\tau_{T}$ and a furrow for $\tau_{q}>\tau_{T}$ respectively
along the trajectory. And both the ridgeline and the furrow become
more and more conspicuous as the difference between $\tau_{q}$ and
$\tau_{T}$ grows. As for the effects of the angular velocity, we have
that a faster angular velocity leads to a lower temperature peak. In
the case of CT with $\tau_{q}>\tau_{T}$, it also leads to the
temperature contours resembling slenderer curved ears of a rabbit.

Apparently, the present work can help us to better understand the
intrinsic mechanism of the DPL heat transfer subjected to a moving
heat source with curved trajectory, and promote the practical
application of the DPL model for the moving heat source problem.
\section*{Acknowledgements}
This work was partially supported by the National Natural Science
Foundation of China, No. 12171240.
The work of the first author was also partially supported by
Postgraduate Research \& Practice Innovation Program of NUAA,
No. xcxjh20220803.

	\bibliographystyle{elsarticle-num} 
	\bibliography{ref/references,ref/localrevised,ref/movingSource}
\end{document}